\documentclass{article}
\usepackage{amssymb}

\usepackage{graphicx}
\usepackage{amsmath}

\newtheorem{theorem}{Theorem}

\newtheorem{corollary}[theorem]{Corollary}

\newtheorem{definition}[theorem]{Definition}

\newtheorem{lemma}[theorem]{Lemma}

\newtheorem{proposition}[theorem]{Proposition}
\newtheorem{remark}[theorem]{Remark}

\input{tcilatex}

\begin{document}

\title{\textbf{Extended-valued topical and anti-topical functions on semimodules}}
\author{Ivan Singer \\
Institute of Mathematics, P. O. Box 1-764, Bucharest, Romania\\
E-mail address: ivan.singer@imar.ro \and Viorel Nitica \\
Department of Mathematics, West Chester University, \\
West Chester, PA 19383 and Institute of Mathematics,\\
P. O. Box 1-764, Bucharest, Romania \\
E-mail address: vnitica@wcupa.edu}
\date{}
\maketitle

\begin{abstract}
In the papers [I. Singer, Lin. Alg. Appl. 2010] and [I. Singer and V.
Nitica, Lin. Alg. Appl. 2012] we have studied functions defined on a \emph{b}%
-complete idempotent semimodule $X$ over a \emph{b}-complete idempotent
semifield $\mathcal{K}=(\mathcal{K},\oplus ,\otimes ),$ with values in $%
\mathcal{K}$, where $\mathcal{K}$ may (or may not) contain a greatest
element $\sup \mathcal{K}$, and the residuation $x/y$ is not defined for $%
x\in X$ and $y=\inf X.$ In the present paper we assume that $\mathcal{K}$
has no greatest element, then adjoin to $\mathcal{K}$\ an outside ``greatest
element'' $\top =\sup \mathcal{K}$ and extend the operations $\oplus $ and $%
\otimes $ from $\mathcal{K}$ to $\overline{\mathcal{K}}:=\mathcal{K}\cup $ $%
\{\top \mathcal{\}}$, so as to obtain a meaning also for $x/\inf X,$ for any 
$x\in X,$ and study functions with values in $\overline{\mathcal{K}}.$ In
fact we consider two different extensions of the product $\otimes $ from $%
\mathcal{K}$ to $\overline{\mathcal{K}}$, denoted by $\otimes $ and $\dot{%
\otimes}$ respectively, and use them to give characterizations of topical
(i.e. increasing homogeneous, defined with the aid of $\otimes $) and
anti-topical (i.e. decreasing anti-homogeneous, defined with the aid of $%
\dot{\otimes}$) functions in terms of some inequalities. Next we introduce
and study for functions $f:X\rightarrow \overline{\mathcal{K}}$ their
conjugates and biconjugates of Fenchel-Moreau type with respect to the
coupling functions $\varphi (x,y)=x/y,$ $\forall x,y\in X,$ and $\psi
(x,(y,d)):=\inf \{x/y,d\},\forall x,y\in X,\forall d\in \overline{\mathcal{K}%
}$, and use them to obtain characterizations of topical and anti-topical
functions. In the subsequent sections we consider for the coupling functions 
$\varphi $ and $\psi $ some concepts that have been studied in [A. M.
Rubinov and I. Singer, Optimization, 2001] and [I. Singer, Optimization,
2004] for the so-called ``additive min-type coupling functions'' $\pi _{\mu
}:R_{\max }^{n}\times R_{\max }^{n}\rightarrow R_{\max }$ and $\pi _{\mu
}:A^{n}\times A^{n}\rightarrow A$ respectively, where $A$ is a conditionally
complete lattice ordered group and $\pi _{\mu }(x,y):=\inf_{1\leq i\leq
n}(x_{i}+y_{i}),\forall x,y\in R_{\max }^{n}$ (or $A^{n}$). Thus, we study
the polars of a set $G\subseteq X$ for the coupling functions $\varphi $ and 
$\psi ,$ and we consider for a function $f:X\rightarrow \overline{\mathcal{K}%
}$ the notion of support set of $f$ with respect to the set $\widetilde{%
\mathcal{T}}$ of all ``elementary topical functions'' $\widetilde{t}%
_{y}(x):=x/y,\forall x\in X,\forall y\in X\backslash \{\inf X\}$ and two
concepts of support set of $f$ at a point $x_{0}\in X.$ The main differences
between the properties of the conjugations with respect to the coupling
functions $\varphi ,\psi $ and $\pi _{\mu }$ and between the properties of
the polars of a set $G$ with respect to the coupling functions $\varphi
,\psi $ and $\pi _{\mu }$ are caused by the fact that while $\pi _{\mu }$ is
symmetric, with values only in $R_{\max }$ (resp. $A),$ $\varphi $ and $\psi 
$ are not symmetric and take values also outside $R_{\max }$ (resp. $A).$

\emph{Key Words}: Semifield, semimodule, \emph{b}-complete, extended
product, extended-valued function, topical function, anti-topical function,
elementary topical function, Fenchel-Moreau conjugate, biconjugate, support
function, polar, bipolar, support set, downward set, subdifferential

\emph{AMS classification:} Primary 06F07, 26B25; Secondary 52A01, 06F20\emph{%
\ }
\end{abstract}

\section{Introduction\label{s1}}

In the previous papers \cite{elem} and \cite{SN}, attempting to contribute
to the construction of a theory of functional analysis and convex analysis
in semimodules over semifields, we have studied topical functions $%
f:X\rightarrow \mathcal{K}$ and related classes of functions, where $X$ is a 
\emph{b}-complete idempotent semimodule over a \emph{b}-complete idempotent
semifield $\mathcal{K}$. We recall that $f:X\rightarrow \mathcal{K}$ is
called \emph{topical }if it is \emph{increasing} (i.e., the relations $%
x^{\prime },x^{\prime \prime }\in X,x^{\prime }\leq x^{\prime \prime }$
imply $f(x^{\prime })\leq f(x^{\prime \prime }),$ where $\leq $ denotes the
canonical order on $\mathcal{K}$, respectively on $X,$ defined by $\lambda
\leq \mu \Leftrightarrow \lambda \oplus \mu =\mu ,\forall \lambda \in 
\mathcal{K},\forall \mu \in \mathcal{K},$ respectively by $x\leq
y\Leftrightarrow x\oplus y=y,\forall x\in X,\forall y\in X),$ and \emph{%
homogeneous} (i.e., $f(\lambda x)=\lambda f(x)$ for all $x\in X,\lambda \in 
\mathcal{K},$ where $\lambda x:=\lambda \otimes x,\lambda f(x)=\lambda
\otimes f(x);$ the fact that we use the same notations for addition $\oplus $
both in $\mathcal{K}$ and in $X$ and for multiplication $\otimes $ both in $%
\mathcal{K}$ and in $\mathcal{K}\times X$ will lead to no confusion). These
definitions will be used also when $\mathcal{K}$ is replaced by $R=((-\infty
,+\infty ),\oplus =\max ,\otimes =+)$ although it is not a semiring, and $X$
is replaced by $R^{n}$. Let us also recall that an idempotent semiring $%
\mathcal{K}$, that is, a semiring with idempotent addition $\oplus $ (i.e.
such that $\lambda \oplus \lambda =\lambda $ for all $\lambda \in \mathcal{K}
$)\ or an idempotent semimodule $X$\ (over an idempotent semiring $\mathcal{K%
}$) is called \emph{b}-\emph{complete}, if it is closed under the sum $%
\oplus $\ of any subset (order-) bounded from above and the multiplication $%
\otimes $\ distributes over such sums.

As in \cite{elem} and \cite{SN}, we shall make the following \emph{basic
assumptions}:

$(A0^{\prime })$ $\mathcal{K}=(\mathcal{K},\oplus ,\otimes )$ is a $b$-\emph{%
complete idempotent semifield }(i.e., a $b$-complete idempotent semiring in
which every $\mu \in \mathcal{K}\backslash \{\varepsilon \}$\ is invertible
for the multiplication $\otimes ,$\ where $\varepsilon $\ denotes the
neutral element of $(\mathcal{K},\oplus )),$ and the supremum of each
(order-) bounded from above subset of $\mathcal{K}$ belongs to $\mathcal{K}$%
; also, $X$ is a $b$-\emph{complete idempotent semimodule} over $\mathcal{K}$%
\emph{. }In the sequel we shall omit the word ``idempotent''; this will lead
to no confusion.

$(A1)$ For all elements $x\in X$ and $y\in X\backslash \{\inf X\}$ the set $%
\{\lambda \in \mathcal{K}|\lambda y\leq x\}$ is (order-) bounded from above,
where $\leq $ denotes the canonical order on $\mathcal{K}$, respectively on $%
X.$

\begin{remark}
\label{rprima}\emph{a) It is easy to see that }an idempotent semifield $%
\mathcal{K}$ has no greatest element $\sup \mathcal{K}$, unless $\mathcal{K}=%
\mathcal{\{}\varepsilon \}$ or $\mathcal{K}=\{\varepsilon ,e\},$ where $%
\varepsilon $ and $e$\ denote the neutral elements of $(\mathcal{K},\oplus )$
and $(\mathcal{K},\otimes )$ respectively \emph{\cite[p. 27, Remark 2.1.2.5]
{Gaubert}}. \emph{Indeed,,\ let }$\mathcal{K}\neq \mathcal{\{}\varepsilon \}$%
\emph{\ and }$\mathcal{K}\neq \{\varepsilon ,e\}$ \emph{and} \emph{assume, a
contrario, that} $\sup \mathcal{K}\in \mathcal{K}.$ \emph{Then} 
\begin{equation*}
\sup \mathcal{K}=(\sup \mathcal{K})\otimes e\leq (\sup \mathcal{K})\otimes
(\sup \mathcal{K})\leq \sup \mathcal{K},
\end{equation*}
\emph{whence }$(\sup \mathcal{K})\otimes (\sup \mathcal{K})=\sup \mathcal{K}$%
\emph{; since }$(\sup \mathcal{K)}^{-1}\neq \varepsilon ,$ \emph{multiplying
both sides of this equality with }$(\sup \mathcal{K)}^{-1}$ \emph{we obtain} 
$\sup \mathcal{K}=e,$\emph{\ so }$\lambda \leq e$\emph{\ for all }$\lambda
\in \mathcal{K}.$\emph{\ Since }$\mathcal{K}$\emph{\ is a semifield,
replacing here each }$\lambda \in \mathcal{K}\backslash \{\varepsilon \}$%
\emph{\ by }$\lambda ^{-1}$\emph{\ it follows that }$\lambda \geq e$\emph{\
and hence }$\lambda =e,$\emph{\ for all }$\lambda \in \mathcal{K}\backslash
\{\varepsilon \}.$\emph{\ Consequently, }$\mathcal{K}=\{\varepsilon ,e\},$ 
\emph{in contradiction with our assumption. The converse statement is
obvious: if }$\mathcal{K}=\{\varepsilon ,e\},$\emph{\ then }$\sup \mathcal{K}%
=e\in \mathcal{K}.$

In the sequel, without any special mention, we shall assume that $\mathcal{K}
$ has no greatest element, \emph{since for} $\mathcal{K}=\mathcal{\{}%
\varepsilon \}$ \emph{the statements are trivial and for }$\mathcal{K}%
=\{\varepsilon ,e\}$\emph{\ the subsequent results remain valid with similar
but simpler proofs.}

\emph{b)} $\mathcal{K}$ \emph{is commutative, by }$(A0^{\prime })$\emph{\
and Iwasawa's theorem (see e.g. \cite{birkh})}.
\end{remark}

An important example of a pair $(X,\mathcal{K})$ satisfying $(A0^{\prime })$
and $(A1)$ is obtained by taking 
\begin{equation}
X=R_{\max }^{n}:=((R\cup \{-\infty \})^{n},\oplus :=\max ,\otimes :=+)
\label{Rnmax}
\end{equation}
(with $\max $ and $+$ understood componentwise), and $\mathcal{K}:=R_{\max
}^{1}.$ The results of \cite{rubsin} on $R^{n}:=(R,\oplus :=\max ,\otimes
:=+)^{n}$ can and will be expressed in the sequel as results on $R_{\max
}^{n}$ replacing $R$ by $R\cup \{-\infty \}$ endowed with the usual
operations $\max $ and $+.$ Also, as has been observed in \cite{rubsin},
many results on $R^{n}$ remain valid, essentially with the same proofs, for $%
X=R_{\max }^{I}$, the set of all bounded vectors $x=(x_{i})_{i\in I}$ where $%
I$ is an arbitrary index set and $x_{i}\in R,\forall i\in I,\sup_{i\in
I}|x_{i}|<+\infty ,$ endowed with the componentwise semimodule operations $%
x^{\prime }\oplus x^{\prime \prime }:=(\max (x_{i}^{\prime },x_{i}^{\prime
\prime }))_{i\in I},\lambda x=(\lambda x_{i})_{i\in I}$ and the
componentwise order relation $x^{\prime }\leq x^{\prime \prime
}\Leftrightarrow x_{i}^{\prime }\leq x_{i}^{\prime \prime },\forall i\in I.$
\ 

One of the main tools in \cite{elem} and \cite{SN} has been residuation. We
recall that by $(A0^{\prime })$ and $(A1)$, for each $y\in X\backslash
\{\inf X\}$ (hence such that the set $\{\lambda \in \mathcal{K}|\lambda
y\leq x\}$ is (order-) bounded from above, by $(A1)$) there exists the \emph{%
residuation operation} $/$ defined by 
\begin{equation}
x/y:=\max \{\lambda \in \mathcal{K}|\lambda y\leq x\},\quad \quad \forall
x\in X,\forall y\in X\backslash \{\inf X\},  \label{resi}
\end{equation}
where $\max $ denotes a supremum that is attained, and it has, among others,
the following properties (see e.g. \cite{CGQ}):

\begin{equation}
(x/y)y\leq x,\quad \quad \forall x\in X,\forall y\in X\backslash \{\inf
X\},\quad \quad  \label{resid5}
\end{equation}
\begin{equation}
y/y=e,\text{\quad \quad }\forall y\in X\backslash \{\inf X\},  \label{resid4}
\end{equation}
\begin{equation}
x/(\mu y)=\mu ^{-1}(x/y),\quad \quad \forall x\in X,\forall \mu \in \mathcal{%
K}\backslash \{\varepsilon \},\forall y\in X\backslash \{\inf X\}.
\label{resid3/2}
\end{equation}

In \cite{elem} and \cite{SN} we have considered only functions $%
f:X\rightarrow \mathcal{K}$, where $(X,\mathcal{K})$ is a pair satisfying $%
(A0^{\prime })$ and $(A1)$. Therefore under the assumption of Remark \ref
{rprima}a) above $x/\inf X,$ i.e. the residuation $x/y$ of (\ref{resi}) for $%
x\in X$ and $y=\inf X,$ is not defined. In the present paper we shall adjoin
to $\mathcal{K}$\ an outside ``greatest element'' $\sup \mathcal{K}$ which
we shall denote by $\top ,$ and extending in a suitable way the operations $%
\oplus $ and $\otimes $ from $\mathcal{K}$ to $\overline{\mathcal{K}}:=%
\mathcal{K}\cup $ $\{\top \mathcal{\}}$, we shall then study functions $%
f:X\rightarrow \overline{\mathcal{K}}:=\mathcal{K}\cup \{\top \},$ that one
may call ``extended-valued functions''; for example, we shall study topical
(i.e. increasing homogeneous) and anti-topical (i.e. decreasing
anti-homogeneous) functions defined on a semimodule $X$ over $\mathcal{K}$
with values in $\overline{\mathcal{K}}:=\mathcal{K}\cup $ $\{\top \mathcal{\}%
}$. Naturally, the extension of the sum operation to sets $M\subseteq 
\overline{\mathcal{K}}$, which we shall denote again by $\oplus ,$ must be $%
\oplus M:=\sup M$ (in $\mathcal{K}$) if $M$ is bounded from above in $%
\mathcal{K}$ and $\oplus M:=\top =\sup \mathcal{K}$ if $M$ is not bounded
from above in $\mathcal{K}$. Also, generalizing the lower addition $%
\underset{\bullet }{+}$ and upper addition $\overset{\bullet }{+}$ on $%
\overline{R},$ of Moreau (see e.g. \cite{moreau}), we shall give two
different extensions of the product $\otimes $ from $\mathcal{K}$ to $%
\overline{\mathcal{K}}$, denoted respectively by $\otimes $ (which will
cause no confusion) and $\dot{\otimes}$, that are dual to each other in a
certain sense, so as to obtain a meaning also for $x/\inf X,$ with any $x\in
X,$ and in particular for $\inf X/\inf X.$ Let us mention that in the
particular case of the pair $(X,\mathcal{K})=(R_{\max }^{n},R_{\max }^{1})$
such extended products and a table of the values of the residuations $x/y$
for all $x,y\in R_{\max }^{1}\cup \{+\infty \mathcal{\}}$ have been given in 
\cite[Table 1]{AGNS}.

First, in Section \ref{s1/2} we shall introduce and study the extended
addition $\oplus $ and the extended products $\otimes $ and $\dot{\otimes}$
in $\overline{\mathcal{K}}:=\mathcal{K}\cup \{\top \}$ and in Section \ref
{s2} we shall use them to give characterizations of topical and anti-topical
functions $f:X\rightarrow \overline{\mathcal{K}}:=\mathcal{K}\cup \{\top \}$
in terms of some inequalities. For the case of topical functions these
characterizations extend some results of \cite{SN}. Since by \cite[Lemma 2.2]
{sin} each one of the products $\otimes $ and $\dot{\otimes}$ determines
uniquely the other, the simple observation that a function $f:X\rightarrow 
\overline{\mathcal{K}}$ is anti-topical if and only if the function $%
h(x):=f(x)^{-1},\forall x\in X,$ is topical, will permit us to deduce
results on anti-topical functions from those of topical functions.

Another domain where topical and anti-topical functions $f:X\rightarrow 
\overline{\mathcal{K}}:=\mathcal{K}\cup \{\top \}$ play an important role is
that of conjugate functions of Fenchel-Moreau type with respect to coupling
functions $\pi :X\times X\rightarrow \overline{\mathcal{K}},$ defined by $%
f^{c(\pi )}(y):=\sup_{x\in X}f(x)^{-1}\pi (x,y),\forall y\in X.$ In Section 
\ref{s3} the extended products $\otimes $ and $\dot{\otimes}$ will permit us
to introduce for functions $f:X\rightarrow \overline{\mathcal{K}}:=\mathcal{K%
}\cup \{\top \}$ their conjugates and biconjugates with respect to the
coupling functions $\varphi (x,y)=x/y,$ $\forall x\in X,\forall y\in X,$ and 
$\psi (x,(y,d)):=\inf \{x/y,d\},\forall x\in X,\forall y\in X,\forall d\in 
\overline{\mathcal{K}}$, and to use them for the study of topical and
anti-topical functions. We shall also consider the ``lower conjugates'' of $%
f $ with respect to these coupling functions, defined with the aid of the
product $\dot{\otimes},$ that are useful for the study of biconjugates. The
main differences between the properties of the conjugations with respect to
the coupling functions $\varphi ,\psi $ and the so-called ``additive
min-type coupling functions'' $\pi _{\mu }:R_{\max }^{n}\times R_{\max
}^{n}\rightarrow R_{\max },$ resp. $\pi _{\mu }:A^{n}\times A^{n}\rightarrow
A,$ where $A$ is a conditionally complete lattice ordered group, studied
previously e.g. in \cite{rubsin}, respectively \cite{radiant}, defined by $%
\pi _{\mu }(x,y):=\inf_{1\leq i\leq n}(x_{i}\otimes y_{i}),\forall
x=(x_{i})\in R_{\max }^{n}$ (resp. $A^{n}),\forall y=(y_{i})\in R_{\max
}^{n} $ (resp. $A^{n}),$ are caused by the fact that while $\pi _{\mu }$ is
''symmetric'' (i.e., $\pi _{\mu }(x,y)=\pi _{\mu }(y,x),\forall x\in
X,\forall y\in X)$ and takes values only in $R_{\max }$ (resp. $A)$, $%
\varphi $ and $\psi $ are not symmetric and take also the value $+\infty $
(resp. $\top )$ ; for example, since $\pi _{\mu }(x,y)$ is topical both in $%
x $ and in $y,$ while $\varphi (x,y)$ is topical as a function of $x$ and
anti-topical as a function of $y$, it follows that while $f^{c(\pi _{\mu
})}:R_{\max }^{n}\rightarrow \overline{R}$ is always a topical function, the
conjugate function $f^{c(\varphi )}:X\rightarrow \overline{\mathcal{K}}$ is
always anti-topical. Note that the notions of ''conditionally complete
lattice ordered group'' and ''\emph{b}-complete idempotent semifield'' are
equivalent (the only difference is the zero which is not difficult to add to
the first notion).

In the subsequent sections we shall consider for the coupling functions $%
\varphi $ and $\psi $ some concepts that have been studied previously for
the additive min-type coupling functions $\pi _{\mu }:R_{\max }^{n}\times
R_{\max }^{n}\rightarrow R_{\max }$ and $\pi _{\mu }:A^{n}\times
A^{n}\rightarrow A,$ in \cite{rubsin} and \cite{radiant} respectively. Thus,
in Section \ref{s5} we shall study the polars of a set $G\subseteq X$ for
the coupling functions $\varphi $ and $\psi ,$ and in Section \ref{s6} we
shall consider the support set of a function $f:X\rightarrow \overline{%
\mathcal{K}}$ with respect to the set $\widetilde{\mathcal{T}}$ of all
``elementary topical functions'' $\widetilde{t}_{y}(x):=x/y,\forall x\in
X,\forall y\in X\backslash \{\inf X\}$ and two concepts of support set of $%
f:X\rightarrow \overline{\mathcal{K}}$ at a point $x_{0}\in X.$ While for
functions $f:X\rightarrow \overline{\mathcal{K}}$ the theory of conjugations 
$f\rightarrow f^{c(\varphi )}$ is of interest, we shall show that for
subsets $G$ of $X$ the theory of polarities $G\rightarrow G^{o(\varphi )}$
permits to obtain some relevant results. Similarly to the case of conjugates
of functions, the main differences between the properties of the polars of a
set $G$ with respect to the coupling functions $\varphi ,\psi $ and the
additive min-type coupling function $\pi _{\mu }:R_{\max }^{n}\times R_{\max
}^{n}\rightarrow R_{\max },$ are caused by the fact that while $\pi _{\mu }$
is symmetric and takes only values in $R_{\max }$, $\varphi $ and $\psi $
are not symmetric and take also the value $+\infty $ (resp. $\top ).$

\section{Extension of $\mathcal{K}$ to $\overline{\mathcal{K}}=\mathcal{K}%
\cup \{\top \}\label{s1/2}$}

\begin{definition}
\label{denlarge}\emph{Let} $\mathcal{K}=(\mathcal{K},\oplus ,\otimes )$ 
\emph{be a }b\emph{-complete semifield that has no greatest element. } \emph{%
We shall adjoin to }$\mathcal{K}$\emph{\ an outside element, which we shall
denote by }$\top ,$\emph{\ and we shall extend the canonical order }$\leq $%
\emph{\ and the addition} $\oplus $ \emph{from }$\mathcal{K}$ \emph{to} 
\emph{an (canonical) order }$\leq $\emph{\ and an addition} $\oplus $ \emph{%
on }$\overline{\mathcal{K}}=\mathcal{K}\cup \{\top \}$ \emph{by} 
\begin{eqnarray}
\varepsilon \leq \alpha \leq \top ,\quad \quad \forall \alpha \in \overline{%
\mathcal{K}},  \label{adjoin} \\
\alpha \oplus \top =\top \oplus \alpha =\top ,\quad \quad \forall \alpha \in 
\overline{\mathcal{K}};  \label{adjoin3}
\end{eqnarray}
\emph{\ } \emph{hence} \emph{the equivalence }$\alpha \leq \beta
\Leftrightarrow \alpha \oplus \beta =\beta $ \emph{remains valid for all }$%
\alpha ,\beta \in \overline{\mathcal{K}}.$ \emph{Furthermore, we shall
extend the multiplication }$\otimes $\emph{\ from }$\mathcal{K}$ \emph{to }$%
\overline{\mathcal{K}}:=\mathcal{K}\cup \{\top \}$ \emph{to two
multiplications} $\dot{\otimes}$ \emph{and} $\otimes $ \emph{by the
following rules: } 
\begin{eqnarray}
\alpha \dot{\otimes}\beta =\alpha \otimes \beta ,\quad \quad \forall \alpha
\in \mathcal{K},\forall \beta \in \mathcal{K},  \label{adjoin21} \\
\alpha \dot{\otimes}\top =\top \dot{\otimes}\alpha =\top ,\quad \quad
\forall \alpha \in \overline{\mathcal{K}},  \label{adjoin22} \\
\alpha \otimes \top =\top \otimes \alpha =\top ,\quad \quad \forall \alpha
\in \overline{\mathcal{K}}\backslash \{\varepsilon \},  \label{adjoin24} \\
\alpha \otimes \varepsilon =\varepsilon \otimes \alpha =\varepsilon ,\quad
\quad \forall \alpha \in \overline{\mathcal{K}}.  \label{adjoin25}
\end{eqnarray}
\emph{We shall often denote the extended product }$\otimes $ $\emph{also}$%
\emph{\ by concatenation, which will cause no confusion.\ }

\emph{For the inverses in }$\overline{\mathcal{K}}$\emph{\ with respect to }$%
\otimes $\emph{\ we shall make the convention } 
\begin{equation}
\varepsilon ^{-1}:=\top ,\quad \top ^{-1}:=\varepsilon ,  \label{adjoin4}
\end{equation}
\emph{whence, by the above, } 
\begin{eqnarray}
\varepsilon ^{-1}\varepsilon =\top \varepsilon =\varepsilon \neq e,\quad
\varepsilon ^{-1}\dot{\otimes}\varepsilon =\top \dot{\otimes}\varepsilon
=\top \neq e,  \label{adjoin43} \\
\top ^{-1}\top =\varepsilon \top =\varepsilon \neq e,\quad \top ^{-1}\dot{%
\otimes}\top =\varepsilon \dot{\otimes}\top =\top \neq e.  \label{adjoin44}
\end{eqnarray}

\emph{We shall call the set }$\overline{\mathcal{K}}=\mathcal{K}\cup \{\top
\}$\emph{\ endowed with the operations }$\oplus ,\otimes $ \emph{and} $\dot{%
\otimes}$\emph{\ the }minimal enlargement of\emph{\ }$\mathcal{K}$\emph{.}%
\newline
\end{definition}

\begin{remark}
\label{rexten}\emph{a) The product }$\dot{\otimes}$ \emph{on }$\overline{%
\mathcal{K}}$\emph{\ is associative, i.e. we have} 
\begin{equation}
(\alpha \dot{\otimes}\beta )\dot{\otimes}\gamma =\alpha \dot{\otimes}(\beta 
\dot{\otimes}\gamma ),\quad \quad \forall \alpha ,\beta ,\gamma \in 
\overline{\mathcal{K}}.  \label{asso}
\end{equation}
\emph{Indeed, if }$\top $\emph{\ occurs as a term in one side of (\ref{asso}%
), then that side must be equal to }$\top $\emph{\ (by (\ref{adjoin22})) and
hence one of the terms of the other side of (\ref{asso}), too, must be equal
to }$\top $\emph{\ (since if }$\lambda ,\mu \in \mathcal{K}$ \emph{then} $%
\lambda \dot{\otimes}\mu =\lambda \otimes \mu \in \mathcal{K}$ \emph{by (\ref
{adjoin21}),} \emph{so} $\lambda \dot{\otimes}\mu \neq \top ).$\emph{\ On
the other hand, if }$\top $\emph{\ does not occur in any one of the terms of
(\ref{asso}), then (\ref{asso}) holds by the usual associativity of }$%
\otimes $\emph{\ on }$\mathcal{K}$.

\emph{b) By the definition of }$\otimes $ \emph{(for} $\alpha \in \mathcal{K}%
)$ \emph{and by (\ref{adjoin24})} \emph{(for} $\alpha =e),$ \emph{we have} 
\begin{equation}
\alpha \otimes e=e\otimes \alpha =\alpha ,\quad \quad \forall \alpha \in 
\overline{\mathcal{K}};  \label{unit}
\end{equation}
\emph{furthermore, by (\ref{adjoin21}) and (\ref{adjoin22}) we have} 
\begin{equation}
\alpha \dot{\otimes}e=e\dot{\otimes}\alpha =\alpha ,\quad \quad \forall
\alpha \in \overline{\mathcal{K}},  \label{unit2}
\end{equation}
\emph{i.e., }$e$ \emph{is the unit element} \emph{of }$\overline{\mathcal{K}}
$\emph{\ for both products }$\otimes $\emph{\ and }$\dot{\otimes}.$

\emph{c) }$\varepsilon ^{-1}$ \emph{and }$\top ^{-1}$\emph{\ are called
``inverses'' only by abuse of language, as shown by (\ref{adjoin43}) and (%
\ref{adjoin44}).}

\emph{d) We shall see that with the above definition, the notions and
results of \cite{elem, SN} on functions }$f:X\rightarrow \mathcal{K},$ \emph{%
where }$X$\emph{\ is a semimodule over }$\mathcal{K},$ \emph{admit
extensions in the above sense to functions }$f:X\rightarrow \overline{%
\mathcal{K}},$ \emph{for the extended product} $\otimes $ \emph{of (\ref
{adjoin24}), (\ref{adjoin25}).\ Therefore in\ the sequel whenever we shall
refer to a result of \cite{SN} or \cite{elem}, we shall understand, without
any special mention, its extension (using the above conventions) to functions%
} $f:X\rightarrow \overline{\mathcal{K}},$ \emph{for the extended product }$%
\otimes $\emph{\ on }$\overline{\mathcal{K}}.$

\emph{e) Note that }a priori\emph{\ the extended products (actions) }$%
\lambda \otimes x$\emph{\ and }$\lambda \dot{\otimes}x,$\emph{\ where }$%
\lambda \in \overline{\mathcal{K}}$\emph{\ and }$x\in X,$\emph{\ need not be
defined, except that by the definition of a semimodule }$X$\emph{\ over }$%
\mathcal{K}$\emph{, we have } 
\begin{equation}
\lambda \inf X=\lambda \otimes \inf X:=\inf X,\quad \quad \forall \lambda
\in \mathcal{K}.  \label{blum1}
\end{equation}
\end{remark}

\begin{definition}
\label{dextenprod}\emph{We extend formula (\ref{blum1}) to }$\lambda =\top $%
\emph{\ by defining} 
\begin{equation}
\top \inf X=\top \otimes \inf X:=\inf X,  \label{blum1/2}
\end{equation}
\emph{and we define} 
\begin{equation}
\lambda \dot{\otimes}\inf X:=\lambda \otimes \inf X=\inf X,\quad \quad
\forall \lambda \in \overline{\mathcal{K}}.  \label{blum4}
\end{equation}
\end{definition}

\begin{remark}
\label{rchestie}\emph{If one would define in any way the products }$\top x,$%
\emph{\ where }$x\in X,$\emph{\ then for each ``homogeneous'' function }$%
f:X\rightarrow \overline{\mathcal{K}},$\emph{\ in the ``extended'' sense }$%
f(\lambda x)=\lambda f(x)$\emph{,}$\forall x\in X,\forall \lambda \in 
\overline{\mathcal{K}},$\emph{\ we would necessarily have} 
\begin{equation*}
f(\top x)=\left\{ 
\begin{array}{l}
\top f(x)=\top \quad \text{\emph{if} }f(x)\neq \varepsilon  \\ 
\top f(x)=\varepsilon \quad \text{\emph{if} }f(x)=\varepsilon ,
\end{array}
\right. 
\end{equation*}
\emph{so }$f(\top x)$\emph{\ could have only two values, namely either }$%
\top $\emph{\ or} $\varepsilon .$
\end{remark}

\begin{definition}
\label{dextenresid}\emph{We define the }extended residuation\emph{\ in any
semimodule }$X$ \emph{over} $\mathcal{K}$ \emph{by taking }$\sup $ \emph{%
instead of }$\max $ \emph{in }$($\emph{\ref{resi}}$),$ \emph{that is}, 
\begin{equation}
x/y:=\sup \{\lambda \in \mathcal{K}|\lambda y\leq x\},\quad \quad \forall
x\in X,\forall y\in X.  \label{extenresid0}
\end{equation}
\end{definition}

\begin{remark}
\label{rsup}\emph{a) In Definition \ref{dextenresid} conditions }$%
(A0^{\prime })$ \emph{and }$(A1)$ \emph{need not be assumed and using (\ref
{extenresid0}) and (\ref{blum1}) we have } 
\begin{equation}
x/\inf X:=\sup \{\lambda \in \mathcal{K}|\lambda \inf X\leq
x\}=\sup_{\lambda \in \mathcal{K}}\lambda =\top ,\quad \quad \forall x\in X;
\label{extenresid}
\end{equation}
\emph{note also that, by (\ref{extenresid0}),} 
\begin{equation}
\inf X/x:=\sup \{\lambda \in \mathcal{K}|\lambda x\leq \inf X\}=\left\{ 
\begin{array}{l}
\varepsilon \quad \text{\emph{if} }x\neq \inf X\text{ } \\ 
\sup_{\lambda \in \mathcal{K}}\lambda =\top \quad \text{\emph{if} }x=\inf X.
\end{array}
\right.   \label{extenresid2}
\end{equation}
\emph{On the other hand, if }$(A0^{\prime })$ \emph{and }$(A1)$\emph{\ hold,
then the }$\sup $ \emph{of (\ref{extenresid0}) is attained in} $\mathcal{K}$ 
\emph{and thus} $x/y$ \emph{of (\ref{extenresid0}) for }$y\neq \inf X$ \emph{%
coincides with }$x/y$ \emph{of (\ref{resi}), so (\ref{resi}), (\ref
{extenresid}) and (\ref{extenresid2}) could be used as an alternative
definition of the extended residuation }$x/y.$ \emph{Note that under }$%
(A0^{\prime })$ \emph{and }$(A1)$ \emph{we need to take }$\sup $\emph{\
instead of }$\max $ \emph{in (\ref{extenresid}), since the set }$\{\lambda
\in \mathcal{K}|\lambda \inf X\leq x\}$\emph{\ is not bounded in }$\mathcal{K%
}$\emph{\ for any }$x\in X$ \emph{and by Remark \ref{rprima}a) }$\top \notin 
\mathcal{K}.$

\emph{b) For }$x=\inf X,$ \emph{from (\ref{extenresid}) and/or (\ref
{extenresid2}) it follows that} 
\begin{equation}
\inf X/\inf X=\top .  \label{extenresid3}
\end{equation}
\end{remark}

We shall use the following extension of formula (\ref{resid3/2}):

\begin{lemma}
\label{lresid}We have 
\begin{equation}
x/(\mu y)=\mu ^{-1}\dot{\otimes}(x/y)\quad \quad \forall x,y\in X,\mu \in 
\mathcal{K}.  \label{extenresid4}
\end{equation}
\end{lemma}

\textbf{Proof}. For $y\in X\backslash \{\inf X\}$ and $\mu \in \mathcal{K}%
\backslash \{\varepsilon \},$ (\ref{extenresid4}) reduces to (\ref{resid3/2}%
).

For $y\in X\backslash \{\inf X\}$ and $\mu =\varepsilon $ we have $%
x/(\varepsilon y)=x/\inf X=\top $ and $\varepsilon ^{-1}\dot{\otimes}%
(x/y)=\top \dot{\otimes}(x/y)=\top .$

For $y=\inf X$ and each $\mu \in \mathcal{K}$ we have $x/(\mu \inf X)=x/\inf
X=\top ,\forall x\in X,$ and $\mu ^{-1}\dot{\otimes}(x/\inf X)=\mu ^{-1}\dot{%
\otimes}\top =\top ,\forall x\in X.\quad \quad \square $

The products $\otimes $ and $\dot{\otimes}$ on $\overline{\mathcal{K}}$ are
closely related, as shown by the following result of \cite[Lemma 2.2]{sin},
for which we give here a more transparent proof:

\begin{theorem}
\label{tunu}The product $\dot{\otimes}$ is uniquely determined by $\otimes $%
, namely we have 
\begin{equation}
\lambda \dot{\otimes}\mu =(\lambda ^{-1}\otimes \mu ^{-1})^{-1},\quad \quad
\forall \lambda \in \overline{\mathcal{K}},\forall \mu \in \overline{%
\mathcal{K}}.  \label{biriz}
\end{equation}

Equivalently, the product $\otimes $ is uniquely determined by $\dot{\otimes}
$, namely we have 
\begin{equation}
\lambda \otimes \mu =(\lambda ^{-1}\dot{\otimes}\mu ^{-1})^{-1},\quad \quad
\forall \lambda \in \overline{\mathcal{K}},\forall \mu \in \overline{%
\mathcal{K}}.  \label{biriz2}
\end{equation}
\end{theorem}

\textbf{Proof}. Let us prove (\ref{biriz}).

Case 1$^{\circ }.$ Both $\lambda \in \mathcal{K}\backslash \{\varepsilon \}$
and $\mu \in \mathcal{K}\backslash \{\varepsilon \}.$ Then $\lambda \dot{%
\otimes}\mu =\lambda \otimes \mu $ and hence, since $\mathcal{K}\backslash
\{\varepsilon \}$ is a commutative semigroup, we have 
\begin{equation*}
(\lambda \dot{\otimes}\mu )^{-1}=(\lambda \otimes \mu )^{-1}=\mu
^{-1}\otimes \lambda ^{-1}=\lambda ^{-1}\otimes \mu ^{-1}.
\end{equation*}

Case 2$^{\circ }.$ $\lambda =\varepsilon $ and $\mu \in \overline{\mathcal{K}%
}.$ Then 
\begin{eqnarray*}
\varepsilon \dot{\otimes}\mu &=&\left\{ 
\begin{array}{l}
\varepsilon \otimes \mu =\varepsilon \quad \text{if }\mu \neq \top \\ 
\varepsilon \dot{\otimes}\top =\top \quad \text{if }\mu =\top ,
\end{array}
\right. \\
(\varepsilon ^{-1}\otimes \mu ^{-1})^{-1} &=&(\top \otimes \mu
^{-1})^{-1}=\left\{ 
\begin{array}{l}
\top ^{-1}=\varepsilon \quad \text{if }\mu ^{-1}\neq \varepsilon \\ 
\varepsilon ^{-1}=\top \text{\quad if }\mu ^{-1}=\varepsilon ,
\end{array}
\right.
\end{eqnarray*}
whence (\ref{biriz}) holds also for $\lambda =\varepsilon ,\mu \in \overline{%
\mathcal{K}}.$

Case 3$^{\circ }.$ $\lambda =\top $ and $\mu \in \overline{\mathcal{K}}.$
Then $\top \dot{\otimes}\mu =\top $ and 
\begin{equation*}
(\top ^{-1}\otimes \mu ^{-1})^{-1}=(\varepsilon \otimes \mu
^{-1})^{-1}=\varepsilon ^{-1}=\top ,
\end{equation*}
so (\ref{biriz}) holds also for $\lambda =\top ,\mu \in \overline{\mathcal{K}%
}.\quad \quad $

Case 4$^{\circ }.$ $\lambda \in \mathcal{K}\backslash \{\varepsilon \}$ and $%
\mu \in \overline{\mathcal{K}}.$ If $\mu =\varepsilon ,$ then we have (\ref
{biriz}) by Case 2$^{\circ }$ with $\lambda $ and $\mu $ interchanged. If $%
\mu =\top ,$ then we have (\ref{biriz}) by Case 3$^{\circ \text{ }}$with $%
\lambda $ and $\mu $ interchanged. Finally, if $\mu \in \mathcal{K}%
\backslash \{\varepsilon \},$ then we are in case 1$^{\circ }.$

This proves (\ref{biriz}). Finally, interchanging in (\ref{biriz}) $\lambda
,\mu $ with $\lambda ^{-1}$ and $\mu ^{-1}$ respectively, we obtain (\ref
{biriz2}).\quad \quad $\square $

\begin{remark}
\label{runu}\emph{In the particular case when }$\mathcal{K}=R_{\max },%
\overline{\mathcal{K}}=\overline{R}=R_{\max }\cup \{+\infty \},$\emph{\
endowed with the lower addition }$\otimes =\underset{\bullet }{+}$\emph{\
and upper addition }$\dot{\otimes}=\overset{\bullet }{+}$\emph{, the rules (%
\ref{adjoin})--(\ref{adjoin25}) are satisfied (by the definitions of }$%
\underset{\bullet }{+}$\ \emph{and} $\overset{\bullet }{+})$ \emph{and
formula (\ref{biriz}) becomes the following well-known observation of Moreau 
\cite[formula (2.1)]{moreau}: } 
\begin{equation}
\lambda \overset{\bullet }{+}\mu =-(-\lambda \underset{\bullet }{+}-\mu
),\quad \quad \forall \lambda \in \overline{R},\forall \mu \in \overline{R}.
\label{biriz3}
\end{equation}
\end{remark}

In the sequel the following properties of equivalence of some inequalities
involving the extended products $\otimes ,\dot{\otimes}$ on $\overline{%
\mathcal{K}}$ will be useful:

\begin{lemma}
\label{lineq}For all $\lambda ,\mu ,\beta \in \overline{\mathcal{K}}:$

\emph{A) }The inequality 
\begin{equation}
\lambda \mu \leq \beta \quad   \label{first-ineq-345}
\end{equation}
is equivalent to 
\begin{equation}
\beta ^{-1}\mu \leq \lambda ^{-1}.\quad   \label{second-ineq-345}
\end{equation}
\ \ \ \ \ \ \ 

\emph{B)} The inequality 
\begin{equation}
\lambda \dot{\otimes}\mu \geq \beta \;\;\;  \label{first-ineq-345-anti}
\end{equation}
\newline
is equivalent to 
\begin{equation}
\beta ^{-1}\dot{\otimes}\mu \geq \lambda ^{-1}.\quad 
\label{second-ineq-345-anti}
\end{equation}
\end{lemma}

\textbf{Proof}. A) The inequalities (\ref{first-ineq-345}) and (\ref
{second-ineq-345}) are equivalent if $\lambda ,\mu ,\beta \in \mathcal{K}%
\setminus \{\varepsilon \}$ (indeed, this follows immediately from the fact
that $\alpha \alpha ^{-1}=e$ for any $\alpha \in \mathcal{K}\setminus
\{\varepsilon \}).$ Thus it remains to consider the cases when one of $%
\lambda ,\mu $ or $\beta $ is $\varepsilon $ or $\top .$

Case (I): $\lambda =\varepsilon .$ Then (\ref{first-ineq-345}) means that $%
\varepsilon =\varepsilon \mu \leq \beta ,$ which is true for all $\mu ,\beta
,$ and (\ref{second-ineq-345}) means that $\beta ^{-1}\mu \leq \lambda
^{-1}=\top ,$ which is also true for all $\mu ,\beta .$ Hence (\ref
{first-ineq-345}) $\Leftrightarrow (\ref{second-ineq-345}).$ $\quad $

Case (IIa): $\lambda =\top $ and $\mu =\varepsilon .$ Then (\ref
{first-ineq-345}) means that $\varepsilon =\top \varepsilon \leq \beta ,$
which is true for all $\beta ,$ and (\ref{second-ineq-345}) means that $%
\varepsilon =\beta ^{-1}\mu \leq \lambda ^{-1},$ which is also true for all $%
\beta .$ Hence (\ref{first-ineq-345}) $\Leftrightarrow (\ref{second-ineq-345}%
).$

Case (IIb): $\lambda =\top $ and $\mu \neq \varepsilon .$ Then (\ref
{first-ineq-345}) means that $\top =\top \mu \leq \beta ,$ which implies
that $\beta =\top $, whence $\beta ^{-1}\mu =\varepsilon \mu \leq \lambda
^{-1},$ so (\ref{first-ineq-345}) implies (\ref{second-ineq-345}). In the
reverse direction, (\ref{second-ineq-345}) means that $\beta ^{-1}\mu \leq
\top ^{-1}=\varepsilon ,$ whence either $\beta ^{-1}=\varepsilon $ or $\mu
=\varepsilon .$ But, if $\beta ^{-1}=\varepsilon ,$ then $\beta =\top ,$
whence $\lambda \mu \leq \beta $ and, on the other hand, if $\mu
=\varepsilon ,$ then $\lambda \mu =\lambda \varepsilon \leq \beta .$ Thus in
either case (\ref{second-ineq-345}) implies $(\ref{first-ineq-345}).$

Case (III): $\beta =\varepsilon $ or $\beta =\top :$ the proof of the
equivalence (\ref{first-ineq-345}) $\Leftrightarrow (\ref{second-ineq-345})$
reduces to the above proofs of the cases $\lambda =\varepsilon $
respectively $\lambda =\top ,$ since (\ref{first-ineq-345}) and (\ref
{second-ineq-345}) are symmetric (by interchanging $\lambda $ and $\beta $
with $\beta ^{-1}$ and $\lambda ^{-1}$ respectively).

\ Case (IV): $\mu =\varepsilon .$ Then (\ref{first-ineq-345}) means that $%
\varepsilon =\lambda \varepsilon \leq \beta ,$ which is true for all $%
\lambda ,\beta ,$ and (\ref{second-ineq-345}) means that $\varepsilon =\beta
^{-1}\mu \leq \lambda ^{-1},$ which is also true for all $\lambda ,\beta .$
Hence (\ref{first-ineq-345}) $\Leftrightarrow (\ref{second-ineq-345}).$

Case (Va): $\mu =\top $ and $\lambda =\varepsilon .$ Then (\ref
{first-ineq-345}) means that $\varepsilon =\varepsilon \top \leq \beta ,$
which is true for all $\beta ,$ and (\ref{second-ineq-345}) means that $%
\beta ^{-1}\mu \leq \top =\lambda ^{-1},$ which is also true for all $\beta
. $ Hence (\ref{first-ineq-345}) $\Leftrightarrow (\ref{second-ineq-345}).$

Case (Vb): $\mu =\top $ and $\lambda \neq \varepsilon .$ Then (\ref
{first-ineq-345}) means that $\top =\lambda \top \leq \beta ,$ which implies
that $\beta =\top ,$ whence $\beta ^{-1}\mu =\varepsilon \mu \leq \lambda
^{-1},$ so (\ref{first-ineq-345}) \ implies (\ref{second-ineq-345}). In the
reverse direction, (\ref{second-ineq-345}) means that $\beta ^{-1}\top
=\beta ^{-1}\mu \leq \lambda ^{-1}<\top ,$ whence $\beta ^{-1}=\varepsilon ,$
so $\beta =\top .$ Therefore $\lambda \mu \leq \top =\beta ,$ and thus (\ref
{second-ineq-345}) implies (\ref{first-ineq-345}).

B) Using Theorem \ref{tunu}, passing to inverses of both sides and applying
part A), and then again taking inverses in both sides and applying Theorem 
\ref{tunu}, we obtain 
\begin{eqnarray*}
\lambda \dot{\otimes}\mu &=&(\lambda ^{-1}\otimes \mu ^{-1})^{-1}\geq \beta
\Leftrightarrow \lambda ^{-1}\otimes \mu ^{-1}\leq \beta ^{-1} \\
&\Leftrightarrow &\beta \otimes \mu ^{-1}\leq \lambda \Leftrightarrow (\beta
\otimes \mu ^{-1})^{-1}\geq \lambda ^{-1}\Leftrightarrow \beta ^{-1}\dot{%
\otimes}\mu \geq \lambda ^{-1}.\quad \quad \square
\end{eqnarray*}

\begin{remark}
\label{rineq}\emph{a) In general, for }$\lambda ,\mu ,\beta \in \overline{%
\mathcal{K}}$\emph{, the inequality (\ref{first-ineq-345}) is not equivalent
to the inequality } 
\begin{equation*}
\mu \leq \lambda ^{-1}\beta .
\end{equation*}
\emph{A counterexample is obtained by taking }$\lambda =\beta =\varepsilon $%
\emph{\ and }$\mu \in \mathcal{K}\setminus \{\varepsilon \}$\emph{. Indeed,
then (\ref{first-ineq-345}) becomes }$\varepsilon \leq \varepsilon $\emph{,
thus it is true, while the second inequality becomes }$\mu \leq \top
\varepsilon =\varepsilon $\emph{, which is false.}

\emph{b) In general, for }$\lambda ,\mu ,\beta \in \overline{\mathcal{K}}$%
\emph{, the inequality (\ref{first-ineq-345-anti}) is not equivalent to the
inequality } 
\begin{equation*}
\mu \geq \lambda ^{-1}\dot{\otimes}\beta .
\end{equation*}
\emph{A counterexample is obtained by taking }$\lambda =\beta =\top $\emph{\
and }$\mu \in \mathcal{K}\setminus \{\varepsilon \}$\emph{. Indeed, then (%
\ref{first-ineq-345-anti}) becomes }$\top \geq \top $\emph{, which is true,
while the second inequality becomes }$\mu \geq \top $\emph{, which is false.}

\emph{c) One cannot replace in Lemma \ref{lineq} the inequalities by the
opposite ones.} \emph{For example }$\varepsilon \varepsilon \geq \varepsilon 
$ \emph{does not imply that} $\varepsilon ^{-1}\varepsilon \geq \varepsilon
^{-1},$ \emph{since} $\varepsilon ^{-1}\varepsilon =\top \varepsilon
=\varepsilon \not\geq \top =\varepsilon ^{-1}.$
\end{remark}

\begin{corollary}
\label{ctipmor}For any $\lambda ,\mu ,\nu \in \overline{\mathcal{K}}$ we
have the equivalence 
\begin{equation}
\mu \otimes \nu \leq \lambda \Leftrightarrow \nu \leq \lambda \dot{\otimes}%
(\mu ^{-1}).  \label{ineq}
\end{equation}
\end{corollary}

\textbf{Proof}. By Theorem \ref{tunu}, given any $\lambda ,\mu ,\nu \in 
\overline{\mathcal{K}},$ for the first inequality of (\ref{ineq}) by Lemma 
\ref{lineq}A) we have the equivalence $\mu \otimes \nu \leq \lambda
\Leftrightarrow \lambda ^{-1}\otimes \mu \leq \nu ^{-1}$ and for the second
inequality of (\ref{ineq}) we have the equivalence $\nu \leq \lambda \dot{%
\otimes}(\mu ^{-1})=(\lambda ^{-1}\otimes \mu )^{-1}\Leftrightarrow \nu
^{-1}\geq \lambda ^{-1}\otimes \mu .\quad \quad \square $

\begin{remark}
\label{rdoi} \emph{In the particular case when }$\mathcal{K}=R_{\max },%
\overline{\mathcal{K}}=\overline{R}=R_{\max }\cup \{+\infty \},$\emph{\
endowed with the lower addition }$\otimes =\underset{\bullet }{+}$\emph{\
and upper addition }$\dot{\otimes}=\overset{\bullet }{+}$\emph{, Corollary 
\ref{ctipmor} reduces to \cite[p.119, Proposition 3c]{moreau} (where the
proof is different).}\quad 
\end{remark}

Finally, concerning ``scalar residuation'', note that one can consider the
semimodule $X=\overline{\mathcal{K}}$ over the semiring $\overline{\mathcal{K%
}}$ (which is almost a semifield, since all elements of $\overline{\mathcal{K%
}}\backslash \{\varepsilon ,\top \}$ are invertible), and by Remark \ref
{rsup}a) one can define residuation in this semimodule (briefly, ``in $%
\overline{\mathcal{K}}$'') by (\ref{extenresid0}). Let us give the following
application of Corollary \ref{ctipmor} which shows that residuation in $%
\overline{\mathcal{K}}$ can be expressed with aid of the upper product and
vice versa, the upper product $\lambda \dot{\otimes}\mu $ can be expressed
with aid of residuation in $\overline{\mathcal{K}}$:

\begin{corollary}
\label{cdoi}We have 
\begin{eqnarray}
\lambda /\mu  &=&\lambda \dot{\otimes}(\mu ^{-1}),\quad \quad \forall
\lambda \in \overline{\mathcal{K}},\forall \mu \in \overline{\mathcal{K}},
\label{resid1} \\
\lambda \dot{\otimes}\mu  &=&\lambda /(\mu ^{-1}),\quad \quad \forall
\lambda \in \overline{\mathcal{K}},\forall \mu \in \overline{\mathcal{K}}.
\label{resid2}
\end{eqnarray}
\end{corollary}

\textbf{Proof}. By the definition of residuation in $\overline{\mathcal{K}}$
and by Corollary \ref{ctipmor}, we obtain 
\begin{eqnarray*}
\lambda /\mu &=&\sup \{\nu \in \mathcal{K}|\mu \otimes \nu \leq \lambda
\}=\sup \{\nu \in \mathcal{K}|\nu \leq \lambda \dot{\otimes}(\mu ^{-1})\} \\
&=&\lambda \dot{\otimes}(\mu ^{-1}),\quad \quad \forall \lambda \in 
\overline{\mathcal{K}},\forall \mu \in \overline{\mathcal{K}},\quad \quad
\end{eqnarray*}
which proves (\ref{resid1}). Finally, interchanging in (\ref{resid1}) $\mu $
and $\mu ^{-1}$ we obtain (\ref{resid2}).\quad \quad $\square $

\begin{remark}
\label{rdoijuma}\emph{a) Alternatively, by (\ref{resid3/2}) applied in }$%
\overline{\mathcal{K}}$ \emph{we have} 
\begin{equation}
\lambda /\mu \nu =\mu ^{-1}\dot{\otimes}\lambda /\nu ,\quad \quad \forall
\lambda \in \overline{\mathcal{K}},\forall \mu \in \overline{\mathcal{K}}%
,\forall \nu \in \overline{\mathcal{K}},  \label{resid3}
\end{equation}
\emph{and taking here }$\nu =e,$\emph{\ we obtain (\ref{resid1}).}

\emph{b) In the particular case when }$\mathcal{K}=R_{\max },\overline{%
\mathcal{K}}=\overline{R}=R_{\max }\cup \{+\infty \},$\emph{\ endowed with
the lower addition }$\otimes =\underset{\bullet }{+}$\emph{\ and upper
addition }$\dot{\otimes}=\overset{\bullet }{+}$\emph{, formula (\ref{resid2}%
) has been proved in \cite[Proposition 2.1]{AGNS}.}
\end{remark}

\section{Characterizations of topical and anti-topical functions $%
f:X\rightarrow \overline{\mathcal{K}}$ in terms of some inequalities \label%
{s2}}

\begin{definition}
\label{dtop-antitop}\emph{Let }$(X,\mathcal{K})$\emph{\ be a pair satisying }%
$(A0^{\prime }),(A1),$ \emph{and let }$\overline{\mathcal{K}}=\mathcal{K}%
\cup \{\top \}$\emph{\ be the minimal enlargement of }$\mathcal{K}$\emph{.} 
\emph{A function }$f:X\rightarrow \overline{\mathcal{K}}$\emph{\ is said to
be}

\emph{a)} increasing \emph{(resp.} decreasing\emph{),} \emph{if }$x^{\prime
},x^{\prime \prime }\in X,x^{\prime }\leq x^{\prime \prime }$ \emph{imply} $%
f(x^{\prime })\leq f(x^{\prime \prime })$ \emph{(resp. }$f(x^{\prime })\geq
f(x^{\prime \prime }));$

\emph{b) }homogeneous \emph{(resp. }anti-homogeneous\emph{), if } 
\begin{equation}
f(\lambda x)=\lambda f(x),\quad \quad \forall x\in X,\forall \lambda \in 
\mathcal{K}  \label{hom1}
\end{equation}
$\;$\emph{(resp. if } 
\begin{equation}
f(\lambda x)=\lambda ^{-1}\dot{\otimes}f(x),\quad \quad \forall x\in
X,\forall \lambda \in \mathcal{K)};  \label{antihom}
\end{equation}

\emph{c) }topical\emph{\ (resp. }anti-topical\emph{), if it is increasing
and homogeneous (resp. decreasing and anti-homogeneous).}
\end{definition}

Similarly to the fact that the two products $\otimes $ and $\dot{\otimes}$
on $\overline{\mathcal{K}}$ are closely related, the topical and
anti-topical functions with values in $\overline{\mathcal{K}}$ are closely
related, too, as shown by the following simple lemma:

\begin{lemma}
\label{lsimple}$f:X\rightarrow \overline{\mathcal{K}}$ is anti-topical if
and only if the function $h:X\rightarrow \overline{\mathcal{K}}$ defined by 
\begin{equation}
h(x):=f(x)^{-1},\quad \quad \forall x\in X,  \label{inv}
\end{equation}
is topical.
\end{lemma}

\textbf{Proof}. Clearly $f$ is decreasing if and only if $h$ is increasing.

Furthermore, let $x\in X,\lambda \in \mathcal{K}$. If $f$ is
anti-homogeneous and $x\in X,\lambda \in \mathcal{K}$, then by Theorem \ref
{tunu} we have 
\begin{equation*}
h(\lambda \otimes x)=f(\lambda \otimes x)^{-1}=(\lambda ^{-1}\dot{\otimes}%
f(x))^{-1}=\lambda \otimes f(x)^{-1}=\lambda \otimes h(x),
\end{equation*}
so $h$ is homogeneous. Conversely, if $h$ is homogeneous, then by Theorem 
\ref{tunu} we have 
\begin{equation*}
f(\lambda \otimes x)=h(\lambda \otimes x)^{-1}=(\lambda \otimes
h(x))^{-1}=\lambda ^{-1}\dot{\otimes}h(x)^{-1}=\lambda ^{-1}\dot{\otimes}%
f(x),
\end{equation*}
so $f$ is anti-homogeneous.$\quad \quad \square $

In the sequel we shall often use Lemma \ref{lsimple} and Theorem \ref{tunu}
above.

\begin{lemma}
\label{lanti}\emph{a) }Let $f:X\rightarrow \overline{\mathcal{K}}$ be a
topical function and for each $y\in X$ let $t_{y}:X\rightarrow \overline{%
\mathcal{K}}$ be the function defined by 
\begin{equation}
t_{y}(x)=f(y)x/y,\quad \quad \forall x\in X.  \label{unu}
\end{equation}

Then we have 
\begin{equation}
t_{y}\leq f,\;t_{y}(y)=f(y).  \label{doiprim}
\end{equation}

Conversely, if $t=t_{y}:X\rightarrow \overline{\mathcal{K}}$ is a function
of the form 
\begin{equation}
t(x)=t_{y}(x)=\alpha x/y,\quad \quad \forall x\in X,  \label{doijuma}
\end{equation}
where $y\in X,\alpha \in \overline{\mathcal{K}},$ satisfying \emph{(\ref
{doiprim}),} then $\alpha =f(y),$ so $t_{y}$ is equal to \emph{(\ref{unu})}.

\emph{b) }Let $f$ $:X\rightarrow \overline{\mathcal{K}}$ be an anti-topical
function. For any $y\in X$ define $q_{y}:X\rightarrow \overline{\mathcal{K}}$
by 
\begin{equation}
q_{y}(x):=(x/y)^{-1}\dot{\otimes}f(y),\quad \quad \forall x\in X.
\label{unuanti}
\end{equation}
Then 
\begin{equation}
q_{y}\geq f,\;q_{y}(y)=f(y).  \label{doianti}
\end{equation}

Conversely, if $q=q_{y}:X\rightarrow \overline{\mathcal{K}}$ is a function
of the form 
\begin{equation}
q(x)=q_{y}(x)=(x/y)^{-1}\dot{\otimes}\alpha ,\quad \quad \forall x\in X,
\label{doijumanti}
\end{equation}
where $y\in X,\alpha \in \overline{\mathcal{K}},$ satisfying \emph{(\ref
{doianti}),} then we have $\alpha =f(y),$ so $q_{y}$ is equal to \emph{(\ref
{unuanti})}.
\end{lemma}

\textbf{Proof}. a) Observe first that \emph{for every homogeneous (and hence
for every topical) function }$f:X\rightarrow \overline{\mathcal{K}}$\emph{\
we have} 
\begin{equation}
f(\inf X)=\varepsilon .  \label{bbs-990}
\end{equation}
Indeed, by (\ref{blum1}), the homogeneity of $f$ and (\ref{adjoin25}) we
have 
\begin{equation*}
f(\inf X)=f(\varepsilon \inf X)=\varepsilon f(\inf X)=\varepsilon .
\end{equation*}

Now let $f$ be topical. If $y\in X\backslash \{\inf X\},$ then since $f$ is
topical, by (\ref{unu}), (\ref{resid5}) and (\ref{resid4}) we have 
\begin{equation}
t_{y}(x)=f(y)x/y=f((x/y)y)\leq f(x),\quad \quad \forall x\in X,\quad \quad
\label{plim1}
\end{equation}
\begin{equation}
t_{y}(y)=f(y)(y/y)=f(y)e=f(y).  \label{plim2}
\end{equation}

On the other hand, if $y=\inf X,$ then by (\ref{unu}), (\ref{bbs-990}), (\ref
{extenresid}), (\ref{adjoin25}) and (\ref{adjoin}), 
\begin{equation}
t_{\inf X}(x)=f(\inf X)(x/\inf X)=\varepsilon \top =\varepsilon \leq
f(x),\quad \quad \forall x\in X,  \label{bum1}
\end{equation}
\begin{equation}
t_{\inf X}(\inf X)=f(\inf X)(\inf X/\inf X)=\varepsilon \top =\varepsilon
=f(\inf X).  \label{bum2}
\end{equation}

Conversely, if $t_{y}$ is a function of the form (\ref{doijuma})\emph{\ }%
satisfying (\ref{doiprim}) and $y\in X\backslash \{\inf X\},$ so $y/y=e,$
then we have $\alpha =\alpha y/y=t_{y}(y)=f(y).$ On the other hand, for $%
y=\inf X,$ formulae (\ref{doijuma}), (\ref{extenresid}), (\ref{doiprim}) and
(\ref{bbs-990}) mean that 
\begin{equation}
t_{\inf X}(x)=\alpha x/\inf X=\alpha \top ,\quad \quad \forall x\in X,
\label{plim3}
\end{equation}
\begin{equation}
t_{\inf X}\leq f,\quad t_{\inf X}(\inf X)=f(\inf X)=\varepsilon .
\label{plim4}
\end{equation}
Hence by (\ref{plim3}) for $x=\inf X$ and the second part of (\ref{plim4}),
we obtain 
\begin{equation*}
\alpha \top =t_{\inf X}(\inf X)=\varepsilon ,
\end{equation*}
which, by (\ref{adjoin24}), (\ref{adjoin25}) and (\ref{bbs-990}), implies
that $\alpha =\varepsilon =f(\inf X).$

b) Observe first that \emph{for every anti-homogeneous (and hence for every
anti-topical) function }$f:X\rightarrow \overline{\mathcal{K}}$\emph{\ we
have} 
\begin{equation}
f(\inf X)=\top .  \label{avem}
\end{equation}
Indeed, by Lemma \ref{lsimple} and part a) above, for every anti-homogeneous
function\emph{\ }$f:X\rightarrow \overline{\mathcal{K}}$ and for $%
h:X\rightarrow \overline{\mathcal{K}}$ of (\ref{inv}) we have 
\begin{equation*}
f(\inf X)^{-1}=h(\inf X)=\varepsilon ,
\end{equation*}
whence (\ref{avem}).

Now let $f$ be anti-topical. Then by Lemma \ref{lsimple} $h:X\rightarrow 
\overline{\mathcal{K}}$ of (\ref{inv}) is topical and hence by Theorem \ref
{tunu} and part a) above, for each $x\in X$ and each $y\in X$ we have 
\begin{eqnarray*}
q_{y}(x) &=&(x/y)^{-1}\dot{\otimes}h(y)^{-1}=((x/y)\otimes h(y))^{-1}\geq
h(x)^{-1}=f(x), \\
q_{y}(y) &=&(y/y)^{-1}\dot{\otimes}h(y)^{-1}=((y/y)\otimes
h(y))^{-1}=h(y)^{-1}=f(y).
\end{eqnarray*}

Conversely, if $q_{y}$ is a function of the form (\ref{doijumanti})
satisfying (\ref{doianti}), that is,\textbf{\ }by Theorem \ref{tunu}, if 
\begin{eqnarray*}
((x/y)\otimes \alpha ^{-1})^{-1} &=&(x/y)^{-1}\dot{\otimes}\alpha
=q_{y}(x)\geq f(x),\quad \quad \forall x\in X, \\
((y/y)\otimes \alpha ^{-1})^{-1} &=&q_{y}(y)=f(y),
\end{eqnarray*}
then passing to inverses, we obtain 
\begin{equation*}
(x/y)\otimes \alpha ^{-1}\leq f(x)^{-1}=h(x)\quad \quad \forall x\in X,
\end{equation*}
\begin{equation*}
(y/y)\otimes \alpha ^{-1}=f(y)^{-1}=h(y).
\end{equation*}
Therefore, since by Lemma \ref{lsimple} $h$ is topical, from part a) (with $%
\alpha $ replaced by $\alpha ^{-1})$ we obtain $h(y)=\alpha ^{-1},$ so $%
f(y)=\alpha .\quad \quad \square $

In \cite[Theorem 5]{SN} we have shown that if $(X,\mathcal{K})$ is a pair
satisfying $(A0^{\prime }),(A1)$, then a function $f:X\rightarrow \mathcal{K}
$ is topical if and only if $f(\inf X)=\varepsilon $ and 
\begin{equation}
f(y)x/y\leq f(x),\quad \quad \forall x\in X,\forall y\in X\backslash \{\inf
X\}  \label{ineg00}
\end{equation}
(where the condition on $y$ was needed in order to be able to define $x/y).$
A similar characterization of topical functions $f:X\rightarrow \mathcal{K}$
in which the inequality $f(y)x/y\leq f(x)$ is replaced by $%
f(y)s_{y,d}(x)\leq f(x),\forall x\in X,\forall y\in X\backslash \{\inf
X\},\forall d\in \mathcal{K}$, where 
\begin{eqnarray}
s_{y,d}(x) &:&=\inf \{x/y,d\}=\inf \{\max \{\lambda \in \mathcal{K}|\lambda
y\leq x\},d\},\quad \quad  \notag \\
\forall x &\in &X,\forall y\in X\backslash \{\inf X\},\forall d\in \mathcal{K%
},  \label{aff1}
\end{eqnarray}
has been given in \cite[Theorem 16]{SN}. Now, with the aid of the
multiplication $\otimes $ on $\overline{\mathcal{K}}$ defined above, we
shall extend these results to functions $f:X\rightarrow \overline{\mathcal{K}%
},$ replacing the conditions $y\in X\backslash \{\inf X\}$ and $d\in 
\mathcal{K}$ by $y\in X$ and $d\in \overline{\mathcal{K}}$ respectively. To
this end, for the case of topical functions we shall extend $s_{y,d}(x)$ to
all $y\in X$ and $d\in \overline{\mathcal{K}}$ by (\ref{aff1}) and 
\begin{eqnarray}
s_{\inf X,d}(x) &:&=\inf \{x/\inf X,d\}=\inf \{\top ,d\}=d,\quad \quad
\forall x\in X,\forall d\in \overline{\mathcal{K}},  \label{aff7} \\
s_{y,\top }(x) &:&=\inf \{x/y,\top \}=x/y,\quad \quad \forall x\in X,\forall
y\in X.  \label{aff5}
\end{eqnarray}

Moreover, we shall also give corresponding characterizations of anti-topical
functions, using the multiplication $\dot{\otimes}$ defined above and the
functions 
\begin{equation}
\overline{s}_{y,d}(x):=\sup \{(x/y)^{-1},d\},\quad \quad \forall x\in
X,\forall y\in X,\forall d\in \overline{\mathcal{K}};  \label{aff1-bis}
\end{equation}
note that in particular for the extreme values $d=\varepsilon $ and $d=\top $
in (\ref{aff1-bis}) we have, respectively, 
\begin{equation}
\overline{s}_{y,\varepsilon }(x)=\sup \{(x/y)^{-1},\varepsilon
\}=(x/y)^{-1},\quad \quad \forall x\in X,\forall y\in X,  \label{aoleu}
\end{equation}
\begin{equation}
\overline{s}_{y,\top }(x)=\sup \{(x/y)^{-1},\top \}=\top ,\quad \quad
\forall x\in X,\forall y\in X.  \label{aoleu2}
\end{equation}

\begin{theorem}
\label{lanti-bis1} Let $(X,\mathcal{K})$ be a pair that satisfies $%
(A0^{\prime }),\;(A1)$.

\emph{a)} For a function $f:X\rightarrow \overline{\mathcal{K}}$ the
following statements are equivalent:

\emph{1}$^{\circ }.$ $f$ is topical.

\emph{2}$^{\circ }.$ We have \emph{(\ref{bbs-990})} and 
\begin{equation}
f(y)x/y\leq f(x),\quad \quad \forall x\in X,\forall y\in X.  \label{bbs-993}
\end{equation}

\emph{3}$^{\circ }.$ We have\emph{\ (\ref{bbs-990}) }and 
\begin{equation}
f(y)s_{y,d}(x)\leq f(x),\quad \quad \forall x\in X,\forall y\in X,\forall
d\in \overline{\mathcal{K}}.  \label{bbs-994}
\end{equation}

\emph{b) }For a function $f:X\rightarrow \overline{\mathcal{K}}$ the
following statements are equivalent:

\emph{1}$^{\circ }.$ $f$ is anti-topical.

\emph{2}$^{\circ }.$ We have \emph{(\ref{avem})} and 
\begin{equation}
f(y)\dot{\otimes}(x/y)^{-1}\geq f(x),\quad \quad \forall x\in X,\forall y\in
X.  \label{bbs-995}
\end{equation}

\emph{3}$^{\circ }.$ We have\emph{\ (\ref{avem}) }and 
\begin{equation}
f(y)\dot{\otimes}\overline{s}_{y,d}(x)\geq f(x),\quad \quad \forall x\in
X,\forall y\in X,\forall d\in \overline{\mathcal{K}}.  \label{bbs-994-anti}
\end{equation}
\end{theorem}

\textbf{Proof}. a) The implication 1$^{\circ }\Rightarrow 2^{\circ }$
follows from Lemma \ref{lanti} a) and its proof.

2$^{\circ }\Rightarrow 3^{\circ }.$ Assume 2$^{\circ }.$ If $x\in X,y\in
X\backslash \{\inf X\},$ then for any $d\in \overline{\mathcal{K}}$ we have,
by (\ref{aff1}) and 2$^{\circ },$ 
\begin{equation*}
f(y)s_{y,d}(x)=f(y)\inf \{x/y,d\}\leq f(y)x/y\leq f(x).
\end{equation*}
If $y=\inf X,$ then by (\ref{bbs-990}) we have 
\begin{equation*}
f(\inf X)s_{\inf X,d}(x)=\varepsilon s_{\inf X,d}(x)=\varepsilon \leq
f(x),\quad \quad \forall x\in X,\forall d\in \overline{\mathcal{K}}.
\end{equation*}

3$^{\circ }\Rightarrow 2^{\circ }.$ Assume 3$^{\circ }$. If $x\in X,y\in
X\backslash \{\inf X\}$, then by 3$^{\circ }$ with any $d\geq x/y$ we obtain 
\begin{equation*}
f(y)x/y=f(y)\inf \{x/y,d\}=f(y)s_{y,d}(x)\leq f(x).\quad \quad
\end{equation*}

Finally, if $x\in X,y=\inf X,$ then $f(\inf X)(x/\inf X)=\varepsilon \top
=\varepsilon \leq f(x)$.$\quad \quad $

2$^{\circ }\Rightarrow 1^{\circ }$. Assume 2$^{\circ }.$ We need to show
that $f$ is increasing and homogeneous.

Assume that $x,y\in X,y\leq x$. Then $e\leq x/y$ and due to (\ref{unit}) and
(\ref{bbs-993}) one has 
\begin{equation*}
f(y)=f(y)e\leq f(y)x/y\leq f(x),
\end{equation*}
so $f$ is increasing.

Assume now that $x\in $ $X\backslash \{\inf X\}$ and $\lambda \in \mathcal{K}%
\backslash \{\varepsilon \}$, so $\lambda x\in X\backslash \{\inf X\},x/x=e$
(by (\ref{resid4})) and $\lambda \lambda ^{-1}=e.$ Then by (\ref{bbs-993})
with $y=\lambda x$ we have $f(\lambda x)x/\lambda x\leq f(x),$ whence using
also (\ref{resid3/2}), 
\begin{equation}
f(\lambda x)=\lambda f(\lambda x)(\lambda ^{-1}x)/x=\lambda f(\lambda
x)x/\lambda x\leq \lambda f(x).  \label{hom1-002}
\end{equation}

On the other hand, for $x=\inf X,\lambda \in \mathcal{K},$ we have, by (\ref
{blum1}), (\ref{bbs-990}) and (\ref{adjoin25}), 
\begin{equation}
f(\lambda \inf X)=f(\inf X)=\varepsilon =\lambda \varepsilon =\lambda f(\inf
X).  \label{hom1bis}
\end{equation}

Moreover, if $x\in X,\lambda =\varepsilon $, then by $\varepsilon x=\inf
X,\forall x\in X,$ and (\ref{bbs-990}) we have 
\begin{equation}
f(\varepsilon x)=f(\inf X)=\varepsilon =\varepsilon f(x).  \label{hom3-002}
\end{equation}

Furthermore, if $y\in X\backslash \{\inf X\},\lambda \in \mathcal{K},$ then
by (\ref{resid4}) and (\ref{bbs-993}) with $x=\lambda y$ we have 
\begin{equation}
\lambda f(y)=f(y)(\lambda y)/y\leq f(\lambda y).  \label{hom2-002}
\end{equation}

On the other hand, for $y=\inf X,\lambda \in \mathcal{K},$ we have 
\begin{equation}
\lambda f(\inf X)=\lambda \varepsilon =\varepsilon =f(\lambda \inf X).
\label{hom2-003}
\end{equation}

From (\ref{hom1-002})--(\ref{hom2-003}) it follows that $f$ is homogeneous.

b)\emph{\ }1$^{\circ }\Leftrightarrow 2^{\circ }.$ By Lemma \ref{lsimple} $f$
is anti-topical if and only if the function $h$ of (\ref{inv}) is topical,
which, by part a), is equivalent to $h(\inf X)=\varepsilon $ and 
\begin{equation*}
h(y)\otimes x/y\leq h(x),\quad \quad \forall x\in X,y\in X,
\end{equation*}
that is, to $f(\inf X)=\top $ and 
\begin{equation*}
f(y)^{-1}\otimes x/y\leq f(x)^{-1},\quad \quad \forall x\in X,\forall y\in X.
\end{equation*}
Then, taking inverses in both sides, we get equivalence of the latter
inequality with 
\begin{equation*}
(f(y)^{-1}\otimes (x/y))^{-1}\geq f(x),\quad \quad \forall x\in X,\forall
y\in X.
\end{equation*}
which, by Theorem \ref{tunu} with $\lambda =f(y),\mu =(x/y)^{-1},$ is
equivalent to (\ref{bbs-995}). $\quad $

1$^{\circ }\Leftrightarrow 3^{\circ }.$ by Lemma \ref{lsimple} $f$ is
anti-topical if and only if the function $h$ of (\ref{inv}) is topical,
which, by part a), is equivalent to $f(\inf X)=\top $ and 
\begin{equation*}
f(y)^{-1}\otimes s_{y,d}(x)\leq f(x)^{-1},\quad \quad \forall x\in X,\forall
y\in X,\forall d\in \overline{\mathcal{K}},
\end{equation*}
that is, taking inverses in both sides, to $f(\inf X)=\top $ and 
\begin{equation*}
(f(y)^{-1}\otimes s_{y,d}(x))^{-1}\geq f(x),\quad \quad \forall x\in
X,\forall y\in X,\forall d\in \overline{\mathcal{K}}.
\end{equation*}
But, by Theorem \ref{tunu} and the definitions of $s$ and $\overline{s},$
for \ any $x\in X,y\in X$ and $d\in \overline{\mathcal{K}}$ we have 
\begin{eqnarray}
(f(y)^{-1}\otimes s_{y,d}(x))^{-1} &=&f(y)\dot{\otimes}s_{y,d}(x)^{-1}=f(y)%
\dot{\otimes}(\inf \{x/y,d\})^{-1}  \notag \\
&=&f(y)\dot{\otimes}\sup \{(x/y)^{-1},d^{-1}\}=f(y)\dot{\otimes}\overline{s}%
_{y,d^{-1}}(x),  \label{mula}
\end{eqnarray}
so the above inequalities are equivalent to (\ref{bbs-994-anti}).\quad \quad 
$\square $

\begin{remark}
\label{rcanbe}\emph{In the statements of Theorem \ref{lanti-bis1} one can
replace, equivalently, the inequalities by equalities. For example in
Theorem \ref{lanti-bis1} a) Statement 2}$^{\circ }$\emph{\ can be replaced,
equivalently, by }

\emph{2}$^{\prime }.$ We have\ \emph{(\ref{bbs-990})} and \emph{\ } 
\begin{equation}
\sup_{y\in X}f(y)x/y=f(x),\quad \quad \forall x\in X.\quad \quad 
\label{egalit}
\end{equation}

\emph{Indeed, for each }$x\in X\backslash \{\inf X\}$ \emph{the} $\sup $%
\emph{\ in (\ref{egalit}) is attained at }$y=x.$ \emph{On the other hand,
for }$x=\inf X$ \emph{we have} $\sup_{y\in X}f(y)(\inf X/y)=\varepsilon
=f(\inf X)\;$\emph{(by (\ref{bbs-990})).} \emph{Thus 2}$^{\circ }\Rightarrow
2^{\prime }.$\emph{\ The reverse implication is obvious.}

\emph{Similarly, in Theorem \ref{lanti-bis1} a) Statement 3}$^{\circ }$\emph{%
\ can be replaced, equivalently, by:}

\emph{3}$^{\prime }.$ We have \emph{(\ref{bbs-990})} and 
\begin{equation}
\sup_{(y,d)\in X\times \overline{\mathcal{K}}}f(y)s_{y,d}(x)=f(x),\quad
\quad \forall x\in X.  \label{expr1/2}
\end{equation}

\emph{Indeed, for each }$x\in X\backslash \{\inf X\}$ \emph{the} $\sup $%
\emph{\ in (\ref{expr1/2}) is attained at }$y=x,d\geq e$ \emph{(since then }$%
f(x)s_{x,d}(x)=f(x)\inf \{x/x,d\}=f(x)e=f(x))$\emph{;} \emph{furthermore,
for }$x=\inf X,y\in X\backslash \{\inf X\},$ \emph{so} $\inf X/y=\varepsilon
,$ \emph{and any} $d\in \overline{\mathcal{K}},$ \emph{we have, by (\ref
{extenresid2}) and (\ref{bbs-990}), } 
\begin{equation*}
\sup_{y\in X}f(y)s_{y,d}(\inf X)=\sup_{y\in X}f(y)\inf \{\inf
X/y,d)=\sup_{y\in X}f(y)\varepsilon =\varepsilon =f(\inf X);
\end{equation*}
\emph{finally, if} $x=y=\inf X,$ \emph{then for any }$d\in \overline{%
\mathcal{K}}$ \emph{we have, by (\ref{bbs-990}),} 
\begin{equation*}
f(\inf X)s_{\inf X,d}(\inf X)=\varepsilon s_{\inf X,d}(\inf X)=\varepsilon
=f(\inf X).
\end{equation*}
\emph{Thus 3}$^{\circ }\Rightarrow 3^{\prime }.$\emph{\ The reverse
implication is obvious. The cases of Theorem \ref{lanti-bis1}b) are similar. 
}
\end{remark}

\begin{corollary}
\label{c0}\emph{a)} If $f$ is topical and for some $y\in X$ we have $%
f(y)=\top $, then for each $x\in X$ either $f(x)=\top $ or $x/y=\varepsilon $%
.$\quad $

\emph{b) }If $f$ is anti-topical and for some $y\in X$ we have $%
f(y)=\varepsilon $, then for each $x\in X$ either $f(x)=\varepsilon $ or $%
(x/y)^{-1}=\top $.

\emph{c)} A function $f:X\rightarrow \overline{\mathcal{K}}$ cannot be
simultaneously topical and anti-topical.
\end{corollary}

\textbf{Proof}. a) If $f$ is topical, $x,y\in X$ and $f(y)=\top $, then by
Theorem \ref{lanti-bis1}a), implication 1$^{\circ }\Rightarrow 2^{\circ },$
we have $\top x/y\leq f(x).$ Hence if $f(x)\neq \top ,$ then by (\ref
{adjoin24}) we obtain $x/y=\varepsilon $.

The proof of part b) is similar, mutatis mutandis.

c) This follows from the fact that the values of $f$ at $\inf X$ are $f(\inf
X)=\varepsilon $ and $f(\inf X)=\top $ for topical and anti-topical
functions respectively, and $\top \neq \varepsilon .\quad \quad \square $

Now we can prove the following further characterizations of anti-topical
functions:

\begin{theorem}
\label{ttrei}For a function $f:X\rightarrow \overline{\mathcal{K}}$ the
following statements are equivalent:

\emph{1}$^{\circ }.$ $f$ is anti-topical.

\emph{2}$^{\circ }.$ We have $f(\inf X)=\top $ and 
\begin{equation}
f(x)\otimes x/y\leq f(y),\quad \quad \forall x\in X,\forall y\in X.
\label{new2}
\end{equation}

\emph{3}$^{\circ }.$ We have $f(\inf X)=\top $ and 
\begin{equation}
f(x)\otimes s_{y,d}(x)\leq f(y),\quad \quad \forall x\in X,\forall y\in
X,d\in \overline{\mathcal{K}}.  \label{new2juma}
\end{equation}

\emph{4}$^{\circ }.$ We have $f(\inf X)=\top $ and 
\begin{equation}
f(y)/(x/y)\geq f(x),\quad \quad \forall x\in X,\forall y\in X.  \label{new3}
\end{equation}
\end{theorem}

\textbf{Proof}. 1$^{\circ }\Leftrightarrow 2^{\circ }.$ By Lemma \ref
{lsimple} we have 1$^{\circ }$ if and only if the function $h$ of (\ref{inv}%
) is topical, which, by Theorem \ref{lanti-bis1}a) is equivalent to $h(\inf
X)=\varepsilon $ and $h(y)\otimes x/y\leq h(x),$ that is, to $f(\inf X)=\top 
$ and $f(y)^{-1}\otimes x/y\leq f(x)^{-1},\forall x\in X,\forall y\in X.$
But, by Lemma \ref{lineq}A), the latter is equivalent to $f(\inf X)=\top $
and (\ref{new2}).

1$^{\circ }\Leftrightarrow 3^{\circ }.$ By Lemma \ref{lsimple} we have 1$%
^{\circ }$ if and only if we have $h(\inf X)=\varepsilon $ and $h(y)\otimes
s_{y,d}(x)\leq h(y),\forall x\in X,\forall y\in X,$ that is, to $f(\inf
X)=\top $ and $f(y)^{-1}\otimes s_{y,d}(x)\leq f(x)^{-1},\forall x\in
X,\forall y\in X.$ But, by Lemma \ref{lineq}A), the latter is equivalent to $%
f(\inf X)=\top $ and (\ref{new2juma}).

1$^{\circ }\Leftrightarrow 4^{\circ }.$ By Theorem \ref{lanti-bis1}b) we
have 1$^{\circ }$ if and only if we have $f(\inf X)=\top $ and (\ref{bbs-995}%
). But, by Corollary \ref{cdoi}, we have $f(y)\dot{\otimes}%
(x/y)^{-1}=f(y)/(x/y),\forall x\in X,\forall y\in X,$ which shows that 1$%
^{\circ }\Leftrightarrow 4^{\circ }.\quad \quad \square $

The following Corollary of Theorem \ref{lanti-bis1} gives characterizations
of the functions $f$ that satisfy the inequalities (\ref{bbs-993})-(\ref
{bbs-994-anti}):

\begin{corollary}
\label{cbun}Let $(X,\mathcal{K})$ be a pair that satisfies $(A0^{\prime
}),\;(A1)$.

\emph{a) }For a function $f:X\rightarrow \overline{\mathcal{K}}$ the
following statements are equivalent:

\emph{1}$^{\circ }.$ We have \emph{(\ref{bbs-993})}.

\emph{2}$^{\circ }.$ We have \emph{(\ref{bbs-994})}.

\emph{3}$^{\circ }.$ Either $f$ is topical or $f\equiv \top .$

\emph{b) }For a function $f:X\rightarrow \overline{\mathcal{K}}$ the
following statements are equivalent:

\emph{1}$^{\circ }.$ We have \emph{(\ref{bbs-995})}.

\emph{2}$^{\circ }.$ We have \emph{(\ref{bbs-994-anti})}.

\emph{3}$^{\circ }.$ Either $f$ is anti-topical or $f\equiv \varepsilon .$
\end{corollary}

\textbf{Proof}. a) 1$^{\circ }\Rightarrow 3^{\circ }.$ Assume that we have 1$%
^{\circ }$ and $f$ is not topical, so $f(\inf X)\neq \varepsilon $ (by
Theorem \ref{lanti-bis1}a)). Then by (\ref{extenresid}) and (\ref{bbs-993})
(applied to $y=\inf X)$ we have 
\begin{equation*}
f(\inf X)\top =f(\inf X)x/\inf X\leq f(x),\quad \quad \forall x\in X,
\end{equation*}
whence, since $f(\inf X)\top =\top $ (by $f(\inf X)\neq \varepsilon ),$ it
follows that $\top \leq f(x),\forall x\in X,$ and hence $f\equiv \top .$

3$^{\circ }\Rightarrow 2^{\circ }.$ Assume 3$^{\circ }.$ If $f$ is topical,
then we have (\ref{bbs-994}) by Theorem \ref{lanti-bis1}a). On the other
hand, if $f\equiv \top ,$ then (\ref{bbs-994}) holds since $\top $ is the
greatest element of $\overline{\mathcal{K}}$.

Finally, the implication 2$^{\circ }\Rightarrow 1^{\circ }$ holds since (\ref
{bbs-994}) for $d=\top $ reduces to (\ref{bbs-993}) (by (\ref{aff5})).

b) 1$^{\circ }\Rightarrow 3^{\circ }.$ Assume that we have 1$^{\circ }$ and $%
f$ is not anti-topical, so $f(\inf X)\neq \top $ (by Theorem \ref{lanti-bis1}%
b)). Then by (\ref{extenresid}), $\top ^{-1}=\varepsilon $ and (\ref{bbs-995}%
) applied to $y=\inf X,$ we have 
\begin{equation*}
f(\inf X)\dot{\otimes}\varepsilon =f(\inf X)\dot{\otimes}(x/\inf X)^{-1}\geq
f(x),\quad \quad \forall x\in X,
\end{equation*}
whence, since $f(\inf X)\dot{\otimes}\varepsilon =\varepsilon $ (by $f(\inf
X)\neq \top ),$ it follows that $\varepsilon \geq f(x),\forall x\in X,$ and
hence $f\equiv \varepsilon .$

3$^{\circ }\Rightarrow 2^{\circ }.$ Assume 3$^{\circ }.$ If $f$ is
anti-topical, then we have (\ref{bbs-995}) by Theorem \ref{lanti-bis1}b). On
the other hand, if $f\equiv \varepsilon ,$ then (\ref{bbs-995}) holds since $%
\varepsilon $ is the smallest element of $\overline{\mathcal{K}}$.

Finally, the implication 2$^{\circ }\Rightarrow 1^{\circ }$ holds since (\ref
{bbs-994-anti}) for $d=\varepsilon $ reduces to (\ref{bbs-995}) (by (\ref
{aoleu})).$\quad \quad \square $

If $\alpha \in \overline{\mathcal{K}}$ and $X$ is any set, we shall use the
same notation $\alpha $ also for the constant function $f(x)\equiv \alpha
,\forall x\in X,$ and we shall write briefly $f\equiv \alpha $ or $f=\alpha
. $

\begin{remark}
\label{rTnutop}\emph{The constant function }$f\equiv \top $ \emph{is not
homogeneous, and hence not topical. Indeed, we have }$\top (\varepsilon
x)=\top ,$ \emph{but} $\varepsilon \top (x)=\varepsilon ,\forall x\in X.$ 
\emph{However, }$f\equiv \top $\emph{\ is anti-topical, since }$\top
(\lambda x)=\top ,\lambda ^{-1}\dot{\otimes}\top (x)=\top ,\forall x\in
X,\forall \lambda \in \mathcal{K}.$\emph{\ Let us also mention that the
constant function }$f\equiv \varepsilon $\emph{\ is topical and hence not
anti-topical.}
\end{remark}

\section{Characterizations of topical and anti-topical functions $%
f:X\rightarrow \overline{\mathcal{K}}$ using conjugates of Fenchel-Moreau
type\label{s3}}

We recall that for two sets $X$ and $Y$ and a ``finite coupling function'' $%
\pi :X\times Y\rightarrow R,$ respectively $\pi :X\times Y\rightarrow A,$
where $A=(A,\oplus ,\otimes )$ is a conditionally complete lattice ordered
group, the Fenchel-Moreau conjugate function (with respect to $\pi )$ of a
function $f:X\rightarrow \overline{R},$ respectively $f:X\rightarrow 
\overline{A},$ where $\overline{A}$ is the canonical enlargement of $A,$ has
been studied in \cite{rubsin, radiant}.

For the next characterizations of topical and anti-topical functions $%
f:X\rightarrow \overline{\mathcal{K}}$ it will be convenient to introduce
the following similar notion of Fenchel-Moreau conjugations:

\begin{definition}
\label{dFM}\emph{Let }$(X,\mathcal{K}),(Y,\mathcal{K})$\emph{\ be two pairs
satisfying }$(A0^{\prime }),(A1),$\emph{\ let }$\overline{\mathcal{K}}:=%
\mathcal{K}\cup \{\top \}$ \emph{be the minimal enlargement of }$\mathcal{K}$%
\emph{, and let }$\pi :X\times Y\rightarrow \overline{\mathcal{K}}$ \emph{be
a function, called ``coupling function''. The} Fenchel-Moreau conjugate\emph{%
\ (}with respect to\emph{\ }$\pi )$\emph{\ of a function }$f:X\rightarrow 
\overline{\mathcal{K}}$\emph{\ is the function }$f^{c(\pi )}:Y\rightarrow 
\overline{\mathcal{K}}$\emph{\ defined by } 
\begin{equation}
f^{c(\pi )}(y):=\sup_{x\in X}f(x)^{-1}\pi (x,y),\quad \quad \forall y\in Y.
\label{FM}
\end{equation}
\end{definition}

\begin{remark}
\label{rFM}\emph{a) For the constant function }$f\equiv \top $ \emph{we have}
\begin{equation}
\top ^{c(\pi )}(y)=\varepsilon ,\quad \quad \forall y\in Y.\quad 
\label{FM0}
\end{equation}
\emph{Indeed, by (\ref{adjoin4}) and (\ref{adjoin25}) we have} 
\begin{equation*}
\top ^{c(\pi )}(y)=\sup_{x\in X}(\top (x)^{-1})\pi (x,y)=\sup_{x\in
X}\varepsilon \pi (x,y)=\varepsilon ,\quad \quad \forall y\in Y.\quad \quad 
\end{equation*}

\emph{b) For the constant function }$f\equiv \varepsilon $ \emph{we have} 
\begin{eqnarray}
\varepsilon ^{c(\pi )}(y)=\sup_{x\in X}\varepsilon (x)^{-1}\pi
(x,y)=\sup_{x\in X}\top \pi (x,y)  \notag \\
=\left\{ 
\begin{array}{l}
\top \quad \text{\emph{if} }\exists x_{0}\in X,\pi (x_{0},y)\neq \varepsilon 
\\ 
\varepsilon \quad \text{\emph{if} }\pi (x,y)=\varepsilon ,\forall x\in X.%
\text{ }
\end{array}
\right.   \label{FM1}
\end{eqnarray}

\emph{c) If }$f$ :$X\rightarrow \overline{\mathcal{K}}$ \emph{and} $x_{0}\in
X$ \emph{are such that} $f(x_{0})=\varepsilon $ \emph{(e.g., if }$f$ \emph{%
is homogeneous, or in particular, topical, and }$x_{0}=\inf X$\emph{) and if 
}$y_{0}\in Y,\pi (x_{0},y_{0})\neq \varepsilon ,$\emph{\ then} 
\begin{equation}
f^{c(\pi )}(y_{0})=\top .\text{\quad }\quad   \label{FM2}
\end{equation}

\emph{Indeed, by (\ref{FM}) and }$\varepsilon ^{-1}=\top $\emph{\ we have} $%
f^{c(\pi )}(y_{0})\geq f(x_{0})^{-1}\pi (x_{0},y_{0})=\top \pi (x_{0},y_{0}),
$\emph{\ whence by }$\pi (x_{0},y_{0})\neq \varepsilon ,$ \emph{(\ref
{adjoin25}) and (\ref{adjoin}) we obtain (\ref{FM2}).}

\emph{d) The abstract (axiomatic) theory of Fenchel-Moreau conjugations and
the theory of kernel representations of lower semi-continuous linear
mappings for idempotent semifield-valued functions (see e.g. \cite{LMS} and
the references therein) have developed in parallel and independently. The
close connections between them have been shown in \cite{sin}.}

\emph{e) In general, one does not always have the equalities} 
\begin{equation}
\sup_{x\in X}f(x)^{-1}\pi (x,y)=\sup_{x\in X}\pi (x,y)/f(x),\quad \quad y\in
X,  \label{FM3}
\end{equation}
\emph{and the second term of (\ref{FM3}) is not suitable to define a
conjugate} \emph{function by } $f^{c(\pi )_{2}}(y):=\sup_{x\in X}\pi
(x,y)/f(x)$\emph{,}$\forall y\in X.$\emph{\ Indeed, \ we shall prove these
statements in Remark \ref{rsecond}b) below. Let us only note here that by
Corollary \ref{cdoi} we have} 
\begin{equation}
\pi (x,y)/f(x)=f(x)^{-1}\dot{\otimes}\pi (x,y),\quad \quad \forall x\in
X,\forall y\in X,  \label{FM4}
\end{equation}
\emph{while in Definition \ref{dFM} we have} $f(x)^{-1}\pi (x,y)$ \emph{%
instead of} $f(x)^{-1}\dot{\otimes}\pi (x,y).$
\end{remark}

Returning to Definition \ref{dFM}, here we shall be interested in
Fenchel-Moreau conjugates first in the case where $(X,\mathcal{K})$ is a
pair satifying $(A0^{\prime }),(A1)$ and $Y:=X,$ with the coupling function $%
\pi =\varphi :X\times X\rightarrow \overline{\mathcal{K}}$ defined by 
\begin{equation}
\varphi (x,y):=x/y=y^{\diamond }(x),\quad \quad \forall x\in X,\forall y\in
X.  \label{coupl6}
\end{equation}
In a different context, related to the study of increasing positively
homogeneous functions $f:C\rightarrow R_{+}$ defined on a cone $C$\ of a
locally convex space $X$\ endowed with the order induced by the closure $%
\overline{C}$\ of $C,$\ with values in $R_{+}=R\cup \{+\infty \},$\ the
coupling function (\ref{coupl6}) has been considered in \cite{DMR}.

Second, we shall be interested in Fenchel-Moreau conjugates for the coupling
function $\pi =\psi :X\times (X\times \overline{\mathcal{K}}\mathcal{)}%
\rightarrow \overline{\mathcal{K}}$\emph{\ }defined by\emph{\ } 
\begin{equation}
\psi (x,(y,d)):=\inf \{x/y,d\}=s_{y,d}(x),\quad \quad \forall x\in X,\forall
y\in X,\forall d\in \overline{\mathcal{K}}.  \label{coupl}
\end{equation}

For $\pi =\varphi $ and $\pi =\psi $ Definition \ref{dFM} leads to:

\begin{definition}
\label{dconjs}\emph{If }$(X,\mathcal{K})$\emph{\ is a pair satisfying }$%
(A0^{\prime }),(A1),$\emph{\ }the Fenchel-Moreau conjugate \emph{of a
function }$f:X\rightarrow \overline{\mathcal{K}}$\emph{\ }

\emph{a) }associated to the coupling function $\varphi $\ of\emph{\ (\ref
{coupl6}), or briefly, the }$\varphi $-conjugate of $f,$ \emph{is the
function }$f^{c(\varphi )}:X\rightarrow \overline{\mathcal{K}}$\emph{\
defined by} 
\begin{equation}
f^{c(\varphi )}(y):=\sup_{x\in X}f(x)^{-1}x/y,\quad \quad \forall y\in X;
\label{conj0}
\end{equation}
\emph{\ }

\emph{b)} associated to the coupling function $\psi $\ of\emph{\ (\ref{coupl}%
), or briefly, the }$\psi $-conjugate of $f,$ \emph{is the function }$%
f^{c(\psi )}:X\times \overline{\mathcal{K}}\rightarrow \overline{\mathcal{K}}
$\emph{\ defined by} 
\begin{equation}
f^{c(\psi )}(y,d)=\sup_{x\in X}f(x)^{-1}s_{y,d}(x),\quad \quad \forall y\in
X,\forall d\in \overline{\mathcal{K}}.  \label{conj}
\end{equation}
\end{definition}

\begin{remark}
\label{rsecond}\emph{a) }For any function $f:X\rightarrow \overline{\mathcal{%
K}},$ the conjugate function $f^{c(\varphi )}:X\rightarrow \overline{%
\mathcal{K}}$ is decreasing \emph{(i.e.} $y_{1}\leq y_{2}\Rightarrow
f^{c(\varphi )}(y_{1})\geq f^{c(\varphi )}(y_{2}))$ and anti-homogeneous 
\emph{(i.e. }$f^{c(\varphi )}(\lambda y)=\lambda ^{-1}\dot{\otimes}%
f^{c(\varphi )}(y)$\emph{\ for all }$y\in X$\emph{\ and }$\lambda \in 
\mathcal{K},$\emph{\ since by Lemma \ref{lresid} we have} 
\begin{eqnarray*}
f^{c(\varphi )}(\lambda y)=\sup_{x\in X}f(x)^{-1}x/\lambda y=\sup_{x\in
X}f(x)^{-1}\lambda ^{-1}\dot{\otimes}x/y \\
=\lambda ^{-1}\dot{\otimes}\sup_{x\in X}f(x)^{-1}x/y=\lambda ^{-1}\dot{%
\otimes}f^{c(\varphi )}(y)),
\end{eqnarray*}
\emph{so }$f^{c(\varphi )}$\ is\emph{\ }anti-topical. \emph{This should be
compared with the situation for the conjugates of }$f$ \emph{with respect to
the so-called ``additive min-type coupling functions'' }$\pi _{\mu }:R_{\max
}^{n}\times R_{\max }^{n}\rightarrow R_{\max \text{ }}$\emph{defined \cite
{rubsin,further} by} 
\begin{equation}
\pi _{\mu }(x,y)=\min_{1\leq i\leq n}(x_{i}+y_{i}),\quad \quad \forall
x=(x_{i})\in R_{\max }^{n},\forall y=(y_{i})\in R_{\max }^{n},  \label{scal0}
\end{equation}
\emph{with} $+$ \emph{denoting the usual addition on }$R_{\max },$ \emph{and
respectively }$\pi _{\mu }:A^{n}\times A^{n}\rightarrow A,$ \emph{where }$A$%
\emph{\ is a conditionally complete lattice ordered group,}\ \emph{defined 
\cite{radiant} by} 
\begin{equation}
\pi _{\mu }(x,y):=\inf_{1\leq i\leq n}(x_{i}\otimes y_{i}),\quad \quad
\forall x=(x_{i})\in A^{n},\forall y=(y_{i})\in A^{n}.  \label{scal}
\end{equation}
\emph{\ For example, in the latter case the conjugate }$f^{c(\pi _{\mu
})}:A^{n}\rightarrow \overline{A}$\emph{\ of a function }$f:A^{n}\rightarrow 
\overline{A}$ \emph{(where} $\overline{A}$\emph{\ is the canonical
enlargement of }$A)$ \emph{is a }$\otimes $-topical\emph{\ function (see 
\cite[Proposition 6.1]{radiant}). Note that for }$X=A^{n}$ \emph{and }$%
\varphi ,\pi _{\mu }$\emph{\ of (\ref{coupl6}), (\ref{scal}) we have} 
\begin{equation}
\varphi (x,y)=x/y=\inf_{1\leq i\leq n}(x_{i}\otimes (y_{i}^{-1}))=\pi _{\mu
}(x,y^{-1}),\quad \quad \forall x\in A^{n},\forall y\in A^{n},  \label{scal2}
\end{equation}
\emph{where} $y^{-1}:=(y_{i}^{-1})_{1\leq i\leq n}$ \emph{(see e.g. 
\cite[Remark 2.4(a)]{AGNS} for} $R_{\max }^{n}$\emph{).}

\emph{b) For the coupling function }$\varphi :X\times X\rightarrow \overline{%
\mathcal{K}}$ \emph{of (\ref{coupl6}) and any }$f:X\rightarrow \overline{%
\mathcal{K}}$ \emph{one can define a function} $f^{c(\varphi
)_{2}}:X\rightarrow \overline{\mathcal{K}}$ \emph{by } 
\begin{equation}
f^{c(\varphi )_{2}}(y):=(x/y)/f(x)=f(x)^{-1}\dot{\otimes}x/y,\quad \quad
\forall y\in X.  \label{cephi2}
\end{equation}

\emph{Then for }$f\equiv \top $ \emph{and }$y=\inf X,$ \emph{by (\ref{conj0}%
), (\ref{bbs-990}), (\ref{extenresid3}) and (\ref{cephi2})} \emph{we have} 
\begin{equation}
\top ^{c(\varphi )}(\inf X)=\sup_{x\in X}\top (x)^{-1}x/\inf X=\sup_{x\in
X}\varepsilon x/\inf X=\varepsilon ,  \label{scal3}
\end{equation}
\begin{equation}
\top ^{c(\varphi )_{2}}(\inf X)=\sup_{x\in X}\top (x)^{-1}\dot{\otimes}%
(x/\inf X)=\varepsilon \dot{\otimes}\top =\top .  \label{scal4}
\end{equation}
\emph{Thus} $\top ^{c(\varphi )}(\inf X)\neq \top ^{c(\varphi )_{2}}(\inf X)$
\emph{and} $\top ^{c(\varphi )_{2}}(\inf X)=\top ,$ \emph{which proves the
statements made in Remark \ref{rFM}e), for the particular coupling function }%
$\pi =\varphi $ \emph{(since usually, extending \cite{ACA}, a conjugation} $%
f\rightarrow f^{c(\pi )}$ \emph{should satisfy} $\top ^{c(\pi )}=\varepsilon
).$ \emph{\ }

\emph{c) By (\ref{aff1}) (extended to all }$y\in X)$\emph{\ and (\ref{aff5}),%
} \emph{for any }$f:X\rightarrow \overline{\mathcal{K}}$ \emph{and the
extreme values }$d=\varepsilon $\emph{\ and }$d=\top $\emph{\ in (\ref{conj}%
) we have, respectively,} 
\begin{align}
f^{c(\psi )}(y,\varepsilon )& =\sup_{x\in X}f(x)^{-1}s_{y,\varepsilon
}(x)=\sup_{x\in X}f(x)^{-1}\varepsilon =\varepsilon ,\quad \quad \forall
y\in X,  \label{extreme} \\
f^{c(\psi )}(y,\top )& =\sup_{x\in X}f(x)^{-1}s_{y,\top }(x)=\sup_{x\in
X}f(x)^{-1}x/y=f^{c(\varphi )}(y),\quad \quad \forall y\in X.
\label{extreme2}
\end{align}
\end{remark}

\begin{theorem}
\label{tconj}Let $(X,\mathcal{K})$ be a pair satisfying assumptions $%
(A0^{\prime })$ and $(A1)$. For a function $f:X\rightarrow \overline{%
\mathcal{K}}$ the following statements are equivalent:

\emph{1}$^{\circ }.$ $f$ is topical.

\emph{2}$^{\circ }.$ We have \emph{(\ref{bbs-990})} and 
\begin{equation}
f^{c(\varphi )}(y)=f(y)^{-1},\quad \quad \forall y\in X.  \label{conj1/2}
\end{equation}

\emph{3}$^{\circ }.$ We have \emph{(\ref{bbs-990})} and 
\begin{equation}
f^{c(\varphi )}(y)\leq f(y)^{-1},\quad \quad \forall y\in X.  \label{conjbim}
\end{equation}

\emph{4}$^{\circ }.$ We have \emph{(\ref{bbs-990})} and 
\begin{equation}
f^{c(\psi )}(y,d)=f(y)^{-1},\quad \quad \forall y\in X,\forall d\in 
\overline{\mathcal{K}}\backslash \{\varepsilon \}.  \label{conj2}
\end{equation}

\emph{5}$^{\circ }.$ We have \emph{(\ref{bbs-990})} and 
\begin{equation}
f^{c(\psi )}(y,d)\leq f(y)^{-1},\quad \quad \forall y\in X,\forall d\in 
\overline{\mathcal{K}}\backslash \{\varepsilon \}.  \label{conjbim2}
\end{equation}
\end{theorem}

\textbf{Proof. }1$^{\circ }\Rightarrow 2^{\circ }.$ If 1$^{\circ }$ holds,
then by Theorem \ref{lanti-bis1}a)\ and Lemma \ref{lineq}A) we have (\ref
{bbs-990}) and $f(x)^{-1}x/y\leq f(y)^{-1},\forall x\in X,\forall y\in X$.
Consequently by (\ref{conj0}) we obtain $f^{c(\varphi )}(y)\leq f(y)^{-1}.$

In the reverse direction, by (\ref{conj}) we have 
\begin{equation*}
f^{c(\varphi )}(y)=\sup_{x\in X}f(x)^{-1}x/y\geq f(y)^{-1}y/y,\quad \quad
y\in Y.
\end{equation*}

If $y\not=\inf X$, then $y/y=e$ and $f(y)^{-1}y/y=f(y)^{-1}$, so we obtain (%
\ref{conj1/2}).

If $y=\inf X,$ then by (\ref{conj0}), (\ref{bbs-990}) and (\ref{extenresid3}%
) we have 
\begin{equation*}
f^{c(\varphi )}(\inf X)=\sup_{x\in X}f(x)^{-1}x/\inf X\geq f(\inf
X)^{-1}\inf X/\inf X=\top \top =\top ,
\end{equation*}
whence $f^{c(\varphi )}(\inf X)=\top =\varepsilon ^{-1}=f(\inf X)^{-1}.$

The implication 2$^{\circ }\Rightarrow 3^{\circ }$ is obvious.

3$^{\circ }\Rightarrow 1^{\circ }.$ If 3$^{\circ }$ holds, then for all $%
x,y\in X$ we have, by (\ref{conj0}) and (\ref{conjbim}), $f(x)^{-1}x/y\leq
f^{c(\varphi )}(y)\leq f(y)^{-1}$, whence $f(y)x/y\leq f(x)$ (by Lemma \ref
{lineq}). Hence, by Theorem \ref{lanti-bis1} a), $f$ is topical.

1$^{\circ }\Rightarrow 4^{\circ }.$ If 1$^{\circ }$ holds, then by Theorem 
\ref{lanti-bis1} a) we have (\ref{bbs-990}). Assume now 1$^{\circ }$ and let 
$x\in X,y\in X\backslash \{\inf X\},d\in \mathcal{K}\backslash \{\varepsilon
\}.$ Then by (\ref{aff1}) (extended to all $y\in X)$, 1$^{\circ }$ and
Theorem \ref{lanti-bis1} a) we have 
\begin{equation}
f(y)s_{y,d}(x)=f(y)\inf \{x/y,d\}\leq f(y)x/y\leq f(x),  \label{civir}
\end{equation}
whence $f(x)^{-1}s_{y,d}(x)\leq f(y)^{-1}$ (by Lemma \ref{lineq}). Hence by (%
\ref{conj}) and since $x\in X$ has been arbitrary, we get 
\begin{equation}
f^{c(\psi )}(y,d)=\sup_{x\in X}\{f(x)^{-1}s_{y,d}(x)\}\leq f(y)^{-1}.
\label{da}
\end{equation}

On the other hand, by (\ref{aff1}), $y\in X\backslash \{\inf X\}$ and (\ref
{resid4}) we have 
\begin{equation}
s_{y,d}(dy)=\inf \{dy/y,d)=d.  \label{ceva}
\end{equation}
Hence, by (\ref{conj}) and (\ref{ceva}), for any $f:X\rightarrow \overline{%
\mathcal{K}}$ (not necessarily topical) and any $y\in X\setminus \{\inf
X\},d\in \mathcal{K}\backslash \{\varepsilon \}\;$we have 
\begin{equation*}
f^{c(\psi )}(y,d)\geq f(dy)^{-1}s_{y,d}(dy)=f(dy)^{-1}d,
\end{equation*}
whence, since by 1$^{\circ }$ and $dd^{-1}=e$ there holds 
\begin{equation*}
f(dy)^{-1}d=(df(y))^{-1}d=dd^{-1}f(y)^{-1}=f(y)^{-1},
\end{equation*}
we obtain the opposite inequality to (\ref{da}), and hence (\ref{conj2}) for 
$y\in X\backslash \{\inf X\},d\in \mathcal{K}\backslash \{\varepsilon \}.$ \ 

Assume now that $y=\inf X$ and $d\in \mathcal{K}\backslash \{\varepsilon \}.$
Then by (\ref{extenresid}), (\ref{bbs-990}) and (\ref{adjoin24}), 
\begin{eqnarray*}
f^{c(\psi )}(\inf X,d) &=&\sup_{x\in X}f(x)^{-1}\inf \{x/\inf X,d\} \\
&=&\sup_{x\in X}f(x)^{-1}\inf \{\top ,d\}=\sup_{x\in X}f(x)^{-1}d \\
&\geq &f(\inf X)^{-1}d=\varepsilon ^{-1}d=\top d=\top ,
\end{eqnarray*}
so we have $f^{c(\psi )}(\inf X,d)=\top =\varepsilon ^{-1}=f(\inf X)^{-1},$
and hence (\ref{conj2}) for all $d\in \mathcal{K}\backslash \{\varepsilon
\}. $

Finally, for $d=\top $ we have, by (\ref{extreme2}), $f^{c(\psi )}(y,\top
)=f^{c(\varphi )}(y),\forall y\in X,$ so the equalities of 3$^{\circ }$ for $%
d=\top $ and 2$^{\circ }$ coincide. Hence by 1$^{\circ }$ and the
implication 1$^{\circ }\Rightarrow 2^{\circ }$ proved above, we have (\ref
{conj2}) also for $d=\top .$

The implication 4$^{\circ }\Rightarrow 5^{\circ }$ is obvious.

5$^{\circ }\Rightarrow 1^{\circ }.$ If 5$^{\circ }$ holds, then for all $%
x,y\in X,d\in \overline{\mathcal{K}}\backslash \{\varepsilon \},$ we have,
by (\ref{conj}) and (\ref{conjbim2}), $f(x)^{-1}s_{y,d}(x)\leq f^{c(\psi
)}(y)\leq f(y)^{-1}$, whence $f(y)s_{y,d}(x)\leq f(x)$ (by Lemma \ref{lineq}%
).\textbf{\ }Furthermore, for $d=\varepsilon $ we have $f(x)^{-1}s_{y,%
\varepsilon }(x)=\varepsilon \leq f(y)^{-1},$ whence $f(y)s_{y,\varepsilon
}(x)\leq f(x),\forall x\in X,\forall y\in X$ (by Lemma \ref{lineq}).
Consequently, by Theorem \ref{lanti-bis1}a), $f$ is topical.$\quad \quad
\square $

\begin{remark}
\label{radaug}\emph{a) One cannot add }$d=\varepsilon $ \emph{to statement 4}%
$^{\circ }$\emph{, since by (\ref{extreme}) we have }$f^{c(\psi
)}(y,\varepsilon )=\varepsilon ,\forall y\in X,$\emph{\ which shows that for
a topical function }$f:X\rightarrow \overline{\mathcal{K}}$ \emph{(hence} $%
f\not\equiv \top $ \emph{by (\ref{bbs-990}))} \emph{one cannot have} $%
f^{c(\psi )}(y,\varepsilon )=f(y)^{-1},\forall y\in X.$

\emph{b) Alternatively, one can also prove the implication 3}$^{\circ
}\Rightarrow 2^{\circ }$\emph{\ of Theorem \ref{tconj} as follows}: \emph{%
Assume that 3}$^{\circ }$\emph{\ holds, that is, } 
\begin{equation}
f^{c(\varphi )}(y)=\sup_{x\in X}f(x)^{-1}x/y\leq f(y)^{-1},\quad \quad
\forall y\in X.  \label{obt}
\end{equation}

\emph{If }$y\neq \inf X,$\emph{\ then the sup in (\ref{obt}) is attained at }%
$x=y,$\emph{\ since }$f(y)^{-1}y/y=f(y)^{-1}.$

\emph{Furthermore, if }$y=\inf X$\emph{\ and there exists }$x_{0}\in X$\emph{%
\ such that }$f(x_{0})\neq \top ,$\emph{\ or equivalently, }$%
f(x_{0})^{-1}\neq \varepsilon ,$\emph{\ then } 
\begin{equation*}
\sup_{x\in X}f(x)^{-1}(x/\inf X)\geq f(x_{0})^{-1}(x_{0}/\inf X)=\top ,
\end{equation*}
\emph{whence by (\ref{bbs-990}) we obtain }$\sup_{x\in X}f(x)^{-1}(x/\inf
X)=\top =\varepsilon ^{-1}=f(\inf X)^{-1}.$

\emph{Finally, if }$f\equiv \top ,$\emph{\ then } 
\begin{equation*}
\sup_{x\in X}\top (x)^{-1}x/y=\varepsilon =\top (y)^{-1},\quad \quad \forall
y\in X.
\end{equation*}
\emph{Thus in all cases we have equality in (\ref{obt}).\ }

\emph{c) Theorem \ref{tconj} above should be compared with \cite[Theorem 5.3]
{rubsin} and \cite[Theorem 6.2]{radiant}; according to the latter a function}
$f:A^{n}\rightarrow \overline{A},$ \emph{where }$A$ \emph{and} $\overline{A}$%
\emph{\ are as in Remark \ref{rsecond}a) above, is ``}$\otimes $\emph{%
-topical'' (i.e., increasing and ``}$\otimes $\emph{-homogeneous'') if and
only if for the coupling function }$\pi _{\mu }$\emph{\ of (\ref{scal}) we
have } 
\begin{equation}
f^{c(\pi _{\mu })}(y)=[f(y^{-1})]^{-1},\quad \quad \forall y\in A^{n}.
\label{coupl44}
\end{equation}
\emph{The reason for this discrepancy between (\ref{coupl44}) and Theorem 
\ref{tconj} is shown by formula (\ref{scal2}) above.\ }

\emph{d) From Corollary \ref{cbun}a) one obtains characterizations of the
functions }$f$\emph{\ that satisfy the equalities (\ref{conj1/2})-(\ref
{conjbim2}), namely, }for a function $f:X\rightarrow \overline{\mathcal{K}}$
the following statements are equivalent:\emph{\ }

\emph{1}$^{\circ }.$ We have \emph{(\ref{conj1/2})}.

\emph{2}$^{\circ }.$ We have \emph{(\ref{conjbim}).}

\emph{3}$^{\circ }.$ We have \emph{(\ref{conj2})}.

\emph{4}$^{\circ }.$ We have \emph{(\ref{conjbim2})}.

\emph{5}$^{\circ }.$ Either $f$ is topical or $f\equiv \top .$

\emph{Indeed, the implication 1}$^{\circ }\Rightarrow 2^{\circ }$\emph{\ is
obvious. Furthermore,} \emph{if we have (\ref{conjbim}), so sup}$_{x\in
X}f(x)^{-1}x/y\leq f(y)^{-1},\forall y\in X,$ \emph{then by Lemma \ref{lineq}%
A) we obtain the inequalities (\ref{bbs-993}), whence by Corollary \ref{cbun}%
a), either }$f$\emph{\ is topical or }$f\equiv \top .$\emph{\ Conversely, if 
}$f$ \emph{is topical, then (\ref{conj1/2}) holds by Theorem \ref{tconj},
while if }$f\equiv \top ,$ \emph{then (\ref{conj1/2}) holds since } 
\begin{equation*}
\emph{sup}_{x\in X}\top (x)^{-1}x/y=\emph{sup}_{x\in X}\varepsilon
x/y=\varepsilon =\top (y)^{-1},\quad \quad \forall y\in X.
\end{equation*}
\emph{Thus 1}$^{\circ }\Rightarrow 2^{\circ }\Rightarrow 5^{\circ
}\Rightarrow 1^{\circ }.$

\emph{Furthermore, the implication 3}$^{\circ }\Rightarrow 4^{\circ }$ \emph{%
is obvious, the implication 4}$^{\circ }\Rightarrow 2^{\circ }$ \emph{holds
by 4}$^{\circ }$ \emph{applied to }$d=\top $ \emph{(see (\ref{extreme2})),} 
\emph{and the implication }$2^{\circ }\Rightarrow 5^{\circ }$\emph{\ was
proved above. Finally, if }$5^{\circ }$\emph{\ holds and }$f$ \emph{is
topical,} \emph{then (\ref{conj2}) holds by Theorem \ref{tconj}, while if }$%
f\equiv \top ,$\emph{\ then (\ref{conj2}) holds since} 
\begin{eqnarray*}
\top ^{c(\psi )}(y,d)=\emph{sup}_{x\in X}\top (x)^{-1}\inf \{x/y,d\}=\sup
\varepsilon \inf \{x/y,d\}=\varepsilon \\
=\top (y)^{-1},\quad \quad \forall y\in X,\forall d\in \overline{\mathcal{K}}%
\backslash \{\varepsilon \}.
\end{eqnarray*}
\emph{Thus 3}$^{\circ }\Rightarrow 4^{\circ }\Rightarrow 2^{\circ
}\Rightarrow 5^{\circ }\Rightarrow 3^{\circ },$ \emph{which completes the
proof of the equivalences 1}$^{\circ }\Leftrightarrow ...\Leftrightarrow
5^{\circ }.$ \emph{\ }
\end{remark}

\begin{lemma}
\label{lpartial}Let $(X,\mathcal{K})$ be a pair satisfying assumptions $%
(A0^{\prime })$ and $(A1)$ and let $\varphi :X\times X\rightarrow \overline{%
\mathcal{K}}$ be the coupling function \emph{(\ref{coupl6})}. Then

\emph{a)} For each $y\in X$ the partial function $\varphi (.,y)$ is topical.

\emph{b) }For each $x\in X$ the partial function $\varphi (x,.)$ is
anti-topical.
\end{lemma}

\textbf{Proof}. a) Let $y\in X.$ Then by the properties of extended
residuation, for each $x^{\prime },x^{\prime \prime },x,y\in X$ and $\lambda
\in \mathcal{K}$ we have 
\begin{eqnarray*}
x^{\prime } &\leq &x^{\prime \prime }\Rightarrow \varphi (x^{\prime
},y)=x^{\prime }/y\leq x^{\prime \prime }/y=\varphi (x^{\prime \prime },y),
\\
x &\in &X,\lambda \in \mathcal{K}\Rightarrow \varphi (\lambda x,y)=(\lambda
x)/y=\lambda (x/y)=\lambda \varphi (x,y).
\end{eqnarray*}

b) Let $x\in X.$ Then by the properties of extended residuation, for each $%
y^{\prime },y^{\prime \prime },y,x\in X$ and $\lambda \in \mathcal{K}$ we
have 
\begin{eqnarray*}
y^{\prime } &\leq &y^{\prime \prime }\Rightarrow \varphi (x,y^{\prime
})=x/y^{\prime }\geq x/y^{\prime \prime }=\varphi (x,y^{\prime \prime }), \\
y &\in &X,\lambda \in \mathcal{K}\Rightarrow \varphi (x,\lambda y)=x/\lambda
y=\lambda ^{-1}\dot{\otimes}x/y=\lambda ^{-1}\dot{\otimes}\varphi
(x,y).\quad \quad \square
\end{eqnarray*}

\begin{remark}
\label{rpartial}\emph{a) We know already part a), since in other words it
says that the function }$y^{\diamond }:X\rightarrow \overline{\mathcal{K}}$ 
\emph{defined by (\ref{coupl6}) is topical (by \cite{elem} when }$y\in $%
\emph{\ }$X\backslash \{\inf X\}$\emph{\ and since} $(\inf X)^{\diamond
}=./\inf X\equiv \top $\emph{\ by (\ref{extenresid})}$).$

\emph{b)} \emph{Lemma \ref{lpartial} should be compared with the fact that
for the coupling function }$\pi _{\mu }:A^{n}\times A^{n}\rightarrow A$ 
\emph{of (\ref{scal}) the partial functions }$\pi _{\mu }(.,y)$\emph{\ and }$%
\pi _{\mu }(x,.)$\emph{\ are topical.}
\end{remark}

\begin{definition}
\label{dlowerconj}\emph{If }$(X,\mathcal{K})$\emph{\ is a pair satisfying }$%
(A0^{\prime }),(A1),$\emph{\ }the Fenchel-Moreau lower conjugate \emph{of a
function }$f:X\rightarrow \overline{\mathcal{K}}$

\emph{a)} associated to the coupling function $\varphi $\ of\emph{\ (\ref
{coupl6}), or briefly, the }$\varphi $-lower conjugate of $f,$ \emph{is the
function }$f^{\theta (\varphi )}:X\rightarrow \overline{\mathcal{K}}$\emph{\
defined by} 
\begin{eqnarray}
f^{\theta (\varphi )}(y):=\inf_{x\in X}f(x)^{-1}\dot{\otimes}(x/y)^{-1} 
\notag \\
=\inf_{x\in X}(f(x)\otimes (x/y)^{-1})^{-1},\quad \forall y\in X;
\label{anti-conj}
\end{eqnarray}

\emph{b)} associated to the coupling function $\psi $\ of\emph{\ (\ref{coupl}%
), or briefly, the }$\psi $-lower conjugate of $f,$ \emph{is the function }$%
f^{\theta (\psi )}:X\rightarrow \overline{\mathcal{K}}$\emph{\ defined by} 
\begin{eqnarray}
f^{\theta (\psi )}(y,d):=\inf_{x\in X}f(x)^{-1}\dot{\otimes}\overline{s}%
_{y,d}(x)  \notag \\
=\inf_{x\in X}(f(x)\otimes s_{y,d^{-1}}(x))^{-1},\quad \forall y\in
X,\forall d\in \overline{\mathcal{K}},  \label{anti-conj-sub}
\end{eqnarray}
\emph{with} $\overline{s}_{y,d}(x)$\emph{\ of (\ref{aff1-bis}).}

\emph{Here the last equalities in (\ref{anti-conj}) and (\ref{anti-conj-sub}%
) follow from Theorem \ref{tunu} and formula (\ref{mula}).}
\end{definition}

\begin{remark}
\label{rsecond-bis}\emph{a) }For any function $f:X\rightarrow \overline{%
\mathcal{K}}$ such that $f\not\equiv \varepsilon $ the lower conjugate
function $f^{\theta (\varphi )}:X\rightarrow \overline{\mathcal{K}}$ is
topical. \emph{Indeed, since }$y_{1}\leq y_{2}$ \emph{implies that} $%
(x/y_{1})^{-1}\leq (x/y_{2})^{-1},f^{\theta (\varphi )}$ \emph{of (\ref
{anti-conj}) is increasing. Furthermore, if }$y\in X\backslash \{\inf
X\},\lambda \in \mathcal{K},$ \emph{then by (\ref{resid3/2}) we have} 
\begin{eqnarray*}
f^{\theta (\varphi )}(\lambda y)=\inf_{x\in X}f(x)^{-1}\dot{\otimes}%
(x/\lambda y)^{-1}=\inf_{x\in X}f(x)^{-1}\dot{\otimes}(\lambda ^{-1}x/y)^{-1}
\\
=\lambda \inf_{x\in X}f(x)^{-1}\dot{\otimes}(x/y)^{-1}=\lambda f^{\theta
(\varphi )}(y);
\end{eqnarray*}
\emph{on the other hand, if} $y=\inf X$ \emph{then since by} $f\not\equiv
\varepsilon $ \emph{there exists} $x_{0}\in X$ \emph{such that} $%
f(x_{0})^{-1}\neq \top ,$ \emph{we have} 
\begin{eqnarray*}
f^{\theta (\varphi )}(\lambda \inf X)=f^{\theta (\varphi )}(\inf
X)=\inf_{x\in X}f(x)^{-1}\dot{\otimes}(x/\inf X)^{-1} \\
\leq f(x_{0})^{-1}\dot{\otimes}(x_{0}/\inf X)^{-1}=f(x_{0})^{-1}\dot{\otimes}%
\varepsilon =\varepsilon ,
\end{eqnarray*}
\emph{whence} $f^{\theta (\varphi )}(\lambda \inf X)=\varepsilon =\lambda
f^{\theta (\varphi )}(\inf X),$ \emph{so} $f^{\theta (\varphi )}$ \emph{of (%
\ref{anti-conj}) is homogeneous, and hence topical.}

\emph{However, note that for the constant function }$\varepsilon ,$\emph{\
that is,} $\varepsilon (x)\equiv \varepsilon ,\forall x\in X,$ \emph{we have}
\begin{equation}
\varepsilon ^{\theta (\varphi )}(y)=\varepsilon ^{\theta (\psi )}(y,d)=\top
,\quad \quad \forall y\in X,\forall d\in \overline{\mathcal{K}};\quad \quad
\label{counterp}
\end{equation}
\emph{indeed, } 
\begin{equation*}
\varepsilon ^{\theta (\varphi )}(y)=\inf_{x\in X}\{\varepsilon ^{-1}\dot{%
\otimes}(x/y)^{-1}\}=\inf_{x\in X}\{\top \dot{\otimes}(x/y)^{-1}\}=\top
,\quad \quad \forall y\in X,
\end{equation*}
\emph{\ and the last equality of (\ref{counterp}) follows similarly.
Consequently, by Remark \ref{rTnutop}, }the lower conjugate function $%
\varepsilon ^{\theta (\varphi )}$ is anti-topical.

\emph{b) By (\ref{aff1-bis}), (\ref{aoleu}) and (\ref{aoleu2}), for any }$%
f:X\rightarrow \overline{\mathcal{K}}$ \emph{and the extreme values }$%
d=\varepsilon $\emph{\ and }$d=\top $\emph{\ in (\ref{anti-conj-sub}) we
have, respectively,} 
\begin{eqnarray}
f^{\theta (\psi )}(y,\varepsilon ) &=&\inf_{x\in X}f(x)^{-1}\dot{\otimes}%
\overline{s}_{y,\varepsilon }(x)=\inf_{x\in X}\{f(x)^{-1}\dot{\otimes}%
(x/y)^{-1}\}  \label{extreme-bis} \\
&=&f^{\theta (\varphi )}(y),\quad \forall y\in X,  \notag \\
f^{\theta (\psi )}(y,\top ) &=&\inf_{x\in X}f(x)^{-1}\dot{\otimes}\overline{s%
}_{y,\top }(x)=\inf_{x\in X}f(x)^{-1}\dot{\otimes}\top  \label{extreme2-bis}
\\
&=&\top ,\quad \forall y\in X.  \notag
\end{eqnarray}
\end{remark}

\begin{theorem}
\label{tconj-antitop}Let $(X,\mathcal{K})$ be a pair satisfying assumptions $%
(A0^{\prime })$ and $(A1)$. For a function $f:X\rightarrow \overline{%
\mathcal{K}}$ the following statements are equivalent:

\emph{1}$^{\circ }.$ $f$ is anti-topical.

\emph{2}$^{\circ }.$ We have \emph{(\ref{avem})} and 
\begin{equation}
f^{\theta (\varphi )}(y)\geq f(y)^{-1},\quad \quad \forall y\in X.
\label{ineg1}
\end{equation}

\emph{3}$^{\circ }.$ We have \emph{(\ref{avem})} and 
\begin{equation}
f^{\theta (\varphi )}(y)=f(y)^{-1},\quad \quad \forall y\in X.
\label{conj1/2-topical}
\end{equation}

\emph{4}$^{\circ }.$ We have \emph{(\ref{avem})} and 
\begin{equation}
f^{\theta (\psi )}(y,d)\geq f(y)^{-1},\quad \quad \forall y\in X,\forall
d\in \mathcal{K}.  \label{ineg2}
\end{equation}

\emph{5}$^{\circ }.$ We have \emph{(\ref{avem})} and 
\begin{equation}
f^{\theta (\psi )}(y,d)=f(y)^{-1},\quad \quad \forall y\in X,\forall d\in 
\mathcal{K}.  \label{conj1/2bis}
\end{equation}
\end{theorem}

\textbf{Proof}. 1$^{\circ }\Leftrightarrow 3^{\circ }.$ If 1$^{\circ }$
holds, then by Theorem \ref{lanti-bis1} b) we have (\ref{avem}).
Furthermore, by (\ref{anti-conj}), Theorem \ref{lanti-bis1} b) and Lemma \ref
{lineq}B) we have (\ref{ineg1}).

Let us prove now the opposite inequalities 
\begin{equation}
f^{\theta (\varphi )}(y)\leq f(y)^{-1},\quad \quad \forall y\in X.
\label{mula6}
\end{equation}

If $y\in X\backslash \{\inf X\},$ then by (\ref{anti-conj}) and (\ref{resid4}%
), we have 
\begin{equation*}
f^{\theta (\varphi )}(y)=\inf_{x\in X}\{f(x)^{-1}\dot{\otimes}%
(x/y)^{-1}\}\leq f(y)^{-1}\dot{\otimes}(y/y)^{-1}=f(y)^{-1}.
\end{equation*}

Furthermore, for $y=\inf X$ we have, by (\ref{anti-conj}) and (\ref{avem}), 
\begin{eqnarray*}
f^{\theta (\varphi )}(\inf X) &=&\inf_{x\in X}f(x)^{-1}\dot{\otimes}(x/\inf
X)^{-1}\leq f(\inf X)^{-1}\dot{\otimes}(\inf X/\inf X)^{-1} \\
&=&\top ^{-1}\dot{\otimes}\top ^{-1}=\varepsilon \dot{\otimes}\varepsilon
=\varepsilon =\top ^{-1}=f(\inf X)^{-1}.\quad \quad
\end{eqnarray*}

The implication 3$^{\circ }\Rightarrow 2^{\circ }$ is obvious.

2$^{\circ }\Rightarrow 1^{\circ }.$ If 2$^{\circ }$ holds, then for all $%
x,y\in X$ we have, by (\ref{anti-conj}) and (\ref{ineg1}), $f(x)^{-1}\dot{%
\otimes}(x/y)^{-1}\geq f^{\theta (\varphi )}(y)\geq f(y)^{-1}.$ Hence $f(y)%
\dot{\otimes}(x/y)^{-1}\geq f(x)$ (by Lemma \ref{lineq}) and therefore, by
Theorem \ref{lanti-bis1} b), $f$ is anti-topical.

1$^{\circ }\Rightarrow 5^{\circ }.$ If 1$^{\circ }$ holds, then by Theorem 
\ref{lanti-bis1} b) we have (\ref{avem}). Assume now 1$^{\circ }$ and let $%
x,y\in X,d\in \mathcal{K}.$ Then by (\ref{aff1-bis}), 1$^{\circ }$ and
Theorem \ref{lanti-bis1} b) we have 
\begin{equation*}
f(y)\dot{\otimes}\overline{s}_{y,d}(x)=f(y)\dot{\otimes}\sup
\{(x/y)^{-1},d\}\geq f(y)\dot{\otimes}(x/y)^{-1}\geq f(x),
\end{equation*}
and thus $f(x)^{-1}\dot{\otimes}\overline{s}_{y,d}(x)\geq f(y)^{-1}$ (by
Lemma \ref{lineq}). Hence by (\ref{anti-conj-sub}), and since $x\in X$ was
arbitrary, we get 
\begin{equation}
f^{\theta (\psi )}(y,d)=\inf_{x\in X}\{f(x)^{-1}\dot{\otimes}\overline{s}%
_{y,d}(x)\}\geq f(y)^{-1}.  \label{da2-bis}
\end{equation}

In the opposite direction, by (\ref{anti-conj-sub}) and (\ref{aff1-bis}),
for any $f:X\rightarrow \overline{\mathcal{K}},$ $y\in X\backslash \{\inf
X\} $ and $d\in \mathcal{K}\backslash \{\varepsilon \}$ we have 
\begin{eqnarray*}
f^{\theta (\psi )}(y,d) &\leq &f(d^{-1}y)^{-1}\dot{\otimes}\overline{s}%
_{y,d}(d^{-1}y) \\
&=&f(d^{-1}y)^{-1}\dot{\otimes}\sup \{(d^{-1}y/y)^{-1},d\}=f(d^{-1}y)^{-1}%
\dot{\otimes}d.
\end{eqnarray*}
Hence, since by 1$^{\circ }$ and $d\in \mathcal{K}\backslash \{\varepsilon
\} $ we have 
\begin{equation*}
f(d^{-1}y)^{-1}\dot{\otimes}d=(d\dot{\otimes}f(y))^{-1}\dot{\otimes}d=(d^{-1}%
\dot{\otimes}f(y)^{-1})\dot{\otimes}d=f(y)^{-1},
\end{equation*}
we obtain the opposite inequality to (\ref{da2-bis}), and hence the equality
(\ref{conj1/2bis}) for $y\in X\backslash \{\inf X\},d\in \mathcal{K}%
\backslash \{\varepsilon \}.$

Assume now that $y=\inf X$ and $d\in \mathcal{K}\backslash \{\varepsilon \}$%
. Then by (\ref{extenresid}), (\ref{avem}), (\ref{adjoin21}) and (\ref
{adjoin25}) we have 
\begin{eqnarray*}
f^{\theta (\psi )}(\inf X,d) &=&\inf_{x\in X}f(x)^{-1}\dot{\otimes}\sup
\{(x/\inf X)^{-1},d\} \\
&=&\inf_{x\in X}f(x)^{-1}\dot{\otimes}\sup \{\varepsilon ,d\}=\inf_{x\in
X}f(x)^{-1}\dot{\otimes}d \\
&\leq &f(\inf X)^{-1}\dot{\otimes}d=\top ^{-1}\dot{\otimes}d=\varepsilon 
\dot{\otimes}d=\varepsilon d=\varepsilon =f(\inf X)^{-1},
\end{eqnarray*}
so we obtain the opposite inequality to (\ref{da2-bis}), and hence the
equality (\ref{conj1/2bis}) for $y=\inf X,d\in \mathcal{K}\backslash
\{\varepsilon \}$. Thus, the equality (\ref{conj1/2bis}) holds for all $y\in
X,d\in \mathcal{K}\backslash \{\varepsilon \}$.

Finally, for $y\in X,d=\varepsilon $ we have, by (\ref{extreme-bis}), 1$%
^{\circ }$ and the implication 1$^{\circ }\Rightarrow 2^{\circ }$ proved
above, 
\begin{equation*}
f^{\theta (\psi )}(y,\varepsilon )=f^{\theta (\varphi )}(y)=f(y)^{-1}.
\end{equation*}

The implication 5$^{\circ }\Rightarrow 4^{\circ }$ is obvious.

4$^{\circ }\Rightarrow 1^{\circ }.$ If 4$^{\circ }$ holds, then by (\ref
{anti-conj-sub}) and Lemma \ref{lineq} b), we have (\ref{avem}) and (\ref
{bbs-994-anti}). Consequently, by Theorem \ref{lanti-bis1} b), $f$ is
anti-topical.\quad \quad $\square $

The following Corollary of Theorem \ref{tconj-antitop} gives
characterizations of the functions $f$ that satisfy (\ref{ineg1})-(\ref
{conj1/2bis}):

\begin{corollary}
\label{cbun3}Let $(X,\mathcal{K})$ be a pair that satisfies $(A0^{\prime
}),\;(A1)$. For a function $f:X\rightarrow \overline{\mathcal{K}}$ the
following statements are equivalent:

\emph{1}$^{\circ }.$ We have \emph{(\ref{ineg1}).}

\emph{2}$^{\circ }.$ We have \emph{(\ref{conj1/2-topical})}.

\emph{3}$^{\circ }.$ We have \emph{(\ref{ineg2}).}

\emph{4}$^{\circ }.$ We have \emph{(\ref{conj1/2bis})}.

\emph{5}$^{\circ }.$ Either $f$ is anti-topical or $f\equiv \varepsilon .$
\end{corollary}

\textbf{Proof}. 1$^{\circ }\Rightarrow 5^{\circ }.$ Assume that we have 1$%
^{\circ }$ and $f$ is not anti-topical, so $f(\inf X)\neq \top $ (by Theorem 
\ref{tconj-antitop}). Then by (\ref{anti-conj}) and (\ref{ineg1}) applied to 
$y=\inf X$ we have 
\begin{equation*}
\inf_{x\in X}\{f(x)^{-1}\dot{\otimes}(x/\inf X)^{-1}\}=f^{\theta (\varphi
)}(\inf X)\geq f(\inf X)^{-1},
\end{equation*}
whence, by (\ref{extenresid}) and since $f(\inf X)^{-1}>\varepsilon ,$ it
follows that 
\begin{equation*}
f(x)^{-1}\dot{\otimes}\varepsilon =f(x)^{-1}\dot{\otimes}(x/\inf
X)^{-1}>\varepsilon ,\quad \quad \forall x\in X.
\end{equation*}
Therefore we must have $f(x)^{-1}=\top ,\forall x\in X,$ so $f\equiv
\varepsilon .$

5$^{\circ }\Rightarrow 4^{\circ }.$ Assume 5$^{\circ }.$ If $f$ is
anti-topical, then we have (\ref{conj1/2bis}) by Theorem \ref{tconj-antitop}%
. On the other hand, if $f\equiv \varepsilon ,$ then (\ref{conj1/2bis})
holds since by (\ref{anti-conj-sub}) we have 
\begin{equation*}
\varepsilon ^{\theta (\psi )}(y,d)=\inf_{x\in X}\{\varepsilon (x)^{-1}\dot{%
\otimes}\overline{s}_{y,d}(x)\}=\top =\varepsilon (y)^{-1},\quad \quad
\forall y\in X,\forall d\in \mathcal{K}.\quad \quad
\end{equation*}

The implication 4$^{\circ }\Rightarrow 2^{\circ }$ holds by (\ref
{extreme-bis}) applied to $d=\top .$ Finally, the equivalences 2$^{\circ
}\Leftrightarrow 1^{\circ }$ and 4$^{\circ }\Leftrightarrow 3^{\circ }$ hold
by (\ref{mula6}).\quad \quad $\square $

Now we shall attempt to apply the second conjugates (i.e., conjugates of
conjugates), or briefly, \emph{biconjugates, }of a function $f:X\rightarrow 
\overline{\mathcal{K}},$ for the study of topical and anti-topical
functions. In the particular case of the coupling function'' $\pi _{\mu
}:R_{\max }^{n}\times R_{\max }^{n}\rightarrow R_{\max }$ of (\ref{scal0}),
in \cite{rubsin} it has been shown (see \cite[Theorem 5.4 and Lemma 5.1]
{rubsin}) that $f$ is topical if and only if $f^{c(\pi _{\mu })c(\pi _{\mu
})}=f$ (see also \cite[Theorem 6.3 and formula (6.32)]{radiant} for an
extension from $R_{\max }^{n}$ to $A^{n}).$ This approach has used the
so-called dual mappings of Moreau (see e.g. \cite{moreau}), which we shall
now try to adapt.

We recall that $\overline{\mathcal{K}}^{X}$ denotes the set of all functions 
$f:X\rightarrow \overline{\mathcal{K}}.$

\begin{definition}
\label{dDual}\emph{Let }$(X,\mathcal{K}),(Y,\mathcal{K})$\emph{\ be two
pairs satisfying }$(A0^{\prime }),(A1).$

\emph{a) For any} \emph{coupling function }$\pi :X\times Y\rightarrow 
\overline{\mathcal{K}}$ \emph{the coupling function }$\overline{\pi }%
:Y\times X\rightarrow \overline{\mathcal{K}}$\emph{\ defined by} 
\begin{equation}
\overline{\pi }(y,x):=\pi (x,y),\quad \quad \forall x\in X,\forall y\in
Y\quad \quad  \label{refl}
\end{equation}
\emph{will be called the} reflexion \emph{of} $\pi .$

\emph{b) The }dual\emph{\ of any mapping }$u:\overline{\mathcal{K}}%
^{X}\rightarrow \overline{\mathcal{K}}^{Y}$ \emph{is the mapping }$u^{\prime
}:\overline{\mathcal{K}}^{Y}\rightarrow \overline{\mathcal{K}}^{X}$\emph{\
defined by } 
\begin{equation}
h^{u^{\prime }}:=\inf_{\substack{ g\in \overline{\mathcal{K}}^{X}  \\ %
g^{u}\leq h}}g,\quad \quad \forall h\in \overline{\mathcal{K}}^{Y},\quad
\label{dual}
\end{equation}
\emph{where we write} $g^{u}$ \emph{and} $h^{u^{\prime }}$\emph{\ instead of 
}$u(g)$\emph{\ and }$u^{\prime }(h)$\emph{\ respectively.}

\emph{c) The }bidual \emph{of any mapping }$u:\overline{\mathcal{K}}%
^{X}\rightarrow \overline{\mathcal{K}}^{Y}$ \emph{is the mapping }$%
f\rightarrow (f^{u})^{u^{\prime }}$\emph{\ of }$\overline{\mathcal{K}}^{X}$%
\emph{\ into }$\overline{\mathcal{K}}^{X}.$
\end{definition}

For the Fenchel-Moreau conjugation $u=c(\pi )$ (see Definition \ref{dFM}) we
have

\begin{lemma}
\label{lL5.1}If $(X,\mathcal{K}),(Y,\mathcal{K})$\ are two pairs satisfying $%
(A0^{\prime }),(A1),$\ and $\pi :X\times Y\rightarrow \overline{\mathcal{K}}$%
\ is a coupling function, then 
\begin{equation}
c(\pi )^{\prime }=c(\overline{\pi }).  \label{anticonj}
\end{equation}
\end{lemma}

\textbf{Proof}. By (\ref{FM}), Lemma \ref{lineq}A) and (\ref{refl}), for any 
$g\in \overline{\mathcal{K}}^{X}$ and $h\in \overline{\mathcal{K}}^{Y}$ we
have the equivalences 
\begin{eqnarray}
g^{c(\pi )}(y) &\leq &h(y),\forall y\in Y  \notag \\
&\Leftrightarrow &g(x)^{-1}\pi (x,y)\leq h(y),\forall x\in X,\forall y\in Y 
\notag \\
&\Leftrightarrow &h(y)^{-1}\pi (x,y)\leq g(x),\forall x\in X,\forall y\in Y 
\notag \\
&\Leftrightarrow &h(y)^{-1}\overline{\pi }(y,x)\leq g(x),\forall x\in
X,\forall y\in Y  \notag \\
&\Leftrightarrow &h^{c(\overline{\pi })}(x)\leq g(x),\forall x\in X,\quad
\quad  \label{equiv}
\end{eqnarray}
whence by (\ref{dual}), 
\begin{equation*}
h^{c(\pi )^{\prime }}(x)=\inf_{\substack{ g\in \overline{\mathcal{K}}^{X} 
\\ g^{c(\pi )}\leq h}}g(x)=\inf_{\substack{ g\in \overline{\mathcal{K}}^{X} 
\\ h^{c(\overline{\pi })}\leq g}}g(x)=h^{c(\overline{\pi })}(x),\quad \quad
\forall x\in X.\quad \quad \square
\end{equation*}

\begin{remark}
\label{rlin}\emph{a) In the particular case where }$X=Y$\emph{\ and }$\pi
=\varphi :X\times X\rightarrow \overline{\mathcal{K}}$\emph{\ is the
coupling function (\ref{coupl6}), }the $\overline{\varphi }$-conjugate
function $f^{c(\overline{\varphi })}$ of any\emph{\ }function\emph{\ }$%
f:X\rightarrow \overline{\mathcal{K}}$ such that $f(\inf X)\neq \top $
(e.g., of any topical function $f$) \emph{satisfies\ } 
\begin{equation*}
f^{c(\overline{\varphi })}(y)=\sup_{x\in X}f(x)^{-1}\overline{\varphi }%
(x,y)=\sup_{x\in X}f(x)^{-1}\varphi (y,x)=\sup_{x\in X}f(x)^{-1}(y/x)
\end{equation*}
\begin{equation*}
\geq f(\inf X)^{-1}(y/\inf X)=\top ,\quad \quad \forall y\in X,
\end{equation*}
\emph{whence} $f^{c(\overline{\varphi })}(y)=\top ,\forall y\in X,$ \emph{so}
$f^{c(\overline{\varphi })}$ is anti-topical, \emph{while for} $f\equiv \top 
$ \emph{we have} 
\begin{equation*}
\top ^{c(\overline{\varphi })}(y)=\sup_{x\in X}\top (x)^{-1}(y/x)=\sup_{x\in
X}\varepsilon y/x=\varepsilon ,\quad \quad \forall y\in X,\text{ }
\end{equation*}
\emph{so }$\top ^{c(\overline{\varphi })}$ is topical. \emph{This should be
compared with the facts that for any function }$f:X\rightarrow \overline{%
\mathcal{K}},f^{c(\varphi )}$\emph{\ is anti-topical (see Remark \ref
{rsecond} a)) and for any function }$f:X\rightarrow \overline{\mathcal{K}}$%
\emph{\ such that} $f\not\equiv \varepsilon ,$ $f^{\theta (\varphi )}$\emph{%
\ is topical, while for }$f\equiv \varepsilon ,$ $\varepsilon ^{\theta
(\varphi )}$\emph{\ is anti-topical (see Remark \ref{rsecond-bis}a)). }

\emph{b) The inequalities occurring in (\ref{equiv}) are equivalent to each
of } 
\begin{equation}
\pi (x,y)\leq h(y)g(x),\quad \overline{\pi }(y,x)\leq g(x)h(y),\quad \quad
\forall x\in X,\forall y\in Y,  \label{FY}
\end{equation}
\emph{which might be called} \emph{``generalized Fenchel-Young
inequalities'',} \emph{because of the particular case of the so-called
``natural coupling function'' }$\pi :R\times R\rightarrow R$ \emph{defined by%
} $\pi (x,y):=xy,\;\forall x\in R,\forall y\in R.$\emph{\ }
\end{remark}

In the particular case where $X=Y$ and $\pi :X\times X\rightarrow \overline{%
\mathcal{K}}$ is a symmetric coupling function, that is, 
\begin{equation}
\pi (x,y)=\pi (y,x),\quad \quad \forall x\in X,\forall y\in X,  \label{symm}
\end{equation}
Lemma \ref{lL5.1} reduces to the following:

\begin{corollary}
\label{cor5.1}If $(X,\mathcal{K})$ is a pair satisfying $(A0^{\prime
}),(A1), $\ and $\pi :X\times X\rightarrow \overline{\mathcal{K}}$\ is a
symmetric coupling function, then $c(\pi )$ is ``self-dual'', that is, 
\begin{equation}
c(\pi )^{\prime }=c(\pi ).  \label{symm2}
\end{equation}
\end{corollary}

\begin{remark}
\label{rpartic}\emph{In particular, for }$\mathcal{K}=R_{\max }$ \emph{and }$%
X$\emph{\ an arbitrary set, Corollary \ref{cor5.1} has been obtained in 
\cite[Lemma 5.1]{rubsin}.}
\end{remark}

For the Fenchel-Moreau biconjugates $f^{c(\varphi )c(\varphi )^{\prime
}}:=(f^{c(\varphi )})^{c(\varphi )^{\prime }}$ with respect to the coupling
function $\varphi $ of (\ref{coupl6}) we obtain

\begin{theorem}
\label{tbiconj}If $(X,\mathcal{K})$\ is a pair satisfying $(A0^{\prime
}),(A1),$ then for every function $f:X\rightarrow \overline{\mathcal{K}}$ we
have 
\begin{equation}
f^{c(\varphi )c(\varphi )^{\prime }}\leq f.\quad \quad \quad \quad
\label{cohull}
\end{equation}

If $f:X\rightarrow \overline{\mathcal{K}}$ is topical, then 
\begin{equation}
f^{c(\varphi )c(\varphi )^{\prime }}=f;  \label{cohull1}
\end{equation}
the converse is not true, since \emph{(\ref{cohull1})} is also satisfied
e.g. for the anti-topical function $f\equiv \top .$
\end{theorem}

\textbf{Proof}. For any function $f:X\rightarrow \overline{\mathcal{K}}$,
applying formula (\ref{dual}) to $h=f^{c(\varphi )}$ and $u=c(\varphi )$ we
obtain 
\begin{equation*}
f^{c(\varphi )c(\varphi )^{\prime }}(x):=\inf_{\substack{ g\in \overline{%
\mathcal{K}}^{X}  \\ g^{c(\varphi )}\leq f^{c(\varphi )}}}g(x)\leq
f(x),\quad \quad \forall x\in X,\quad \quad \quad \quad
\end{equation*}
which proves (\ref{cohull}).

Furthermore, observe that for any function $f:X\rightarrow \overline{%
\mathcal{K}}$ we have, using (\ref{anticonj}) with $\pi =\varphi $ of (\ref
{coupl6}), 
\begin{equation}
f^{c(\varphi )c(\varphi )^{\prime }}(x)=(f^{c(\varphi )})^{c(\overline{%
\varphi })}(x),\quad \quad \forall x\in X.  \label{biconj2}
\end{equation}
Now let $f:X\rightarrow \overline{\mathcal{K}}$ be a topical function and
let $h=f^{c(\varphi )}$ and $x\in X.$ Then by (\ref{refl}) for $\pi =\varphi 
$ of (\ref{coupl6}) and by $f^{c(\varphi )}(y)=f(y)^{-1},\forall y\in X$
(see Theorem \ref{tconj}), we have 
\begin{eqnarray}
f^{c(\varphi )c(\overline{\varphi })}(x) &=&h^{c(\overline{\varphi }%
)}(x)=\sup_{y\in X}h(y)^{-1}\overline{\varphi }(y,x)=\sup_{y\in
X}h(y)^{-1}(y/x)  \notag \\
&=&\sup_{y\in X}(f^{c(\varphi )}(y)^{-1}(y/x))=\sup_{y\in X}f(y)(y/x)\geq
f(x)(x/x).\quad \quad  \label{star}
\end{eqnarray}
If $x\neq \inf X,$ then $f(x)x/x=f(x)$, and hence by (\ref{star}), 
\begin{equation*}
f^{c(\varphi )c(\overline{\varphi })}(x)\geq f(x)x/x=f(x).
\end{equation*}
On the other hand, if $x=\inf X,$ then $f(\inf X)=\varepsilon $ (since $f$
is topical), and hence by (\ref{star}), 
\begin{equation*}
f^{c(\varphi )c(\overline{\varphi })}(\inf X)\geq f(\inf X)(\inf X/\inf
X)=f(\inf X)\top =\varepsilon =f(\inf X).
\end{equation*}

Thus for any topical function $f:X\rightarrow \overline{\mathcal{K}}$ we
have $f^{c(\varphi )c(\overline{\varphi })}(x)\geq f(x),$ $\forall x\in X,$
which together with (\ref{cohull}) yields the equality (\ref{cohull1}).

Finally, let us show that the converse implication is not true. Indeed, by
Remarks \ref{rFM}a) and \ref{rlin}a) we have 
\begin{equation*}
\top ^{c(\varphi )c(\overline{\varphi })}=(\top ^{c(\varphi )})^{c(\overline{%
\varphi })}=\varepsilon ^{c(\overline{\varphi })}=\top ,
\end{equation*}
whence (\ref{cohull1}) for $f\equiv \top ,$ which, since $f\equiv \top $ is
anti-topical (see Remark \ref{rTnutop}), completes the proof.\quad \quad $%
\square $

\begin{remark}
\label{rbiconj}\emph{In \cite[Theorem 5.4]{rubsin} it has been shown that
for any topical function }$f:R_{\max }^{n}\rightarrow \overline{R}$ \emph{%
the equality (\ref{cohull1}) holds for }$\varphi $\emph{\ replaced by the
additive min-type coupling function }$\pi _{\mu }$\emph{\ of (\ref{scal0})
and conversely, every function }$f:R_{\max }^{n}\rightarrow \overline{R}$ 
\emph{such that} $f^{c(\pi _{\mu })c(\pi _{\mu })^{\prime }}=f$ \emph{is
topical.} \emph{The proof of \cite{rubsin} has used the fact that }$\pi
_{\mu }$ \emph{is symmetric and hence self-dual, but this is not the case
for }$\varphi ;$ \emph{furthermore,} $f^{c(\pi _{\mu })}$ \emph{is topical
for every function }$f:R_{\max }^{n}\rightarrow \overline{R}$\emph{, but }$%
f^{c(\varphi )}$ \emph{is anti-topical for every function }$f:X\rightarrow 
\overline{\mathcal{K}}$\emph{\ (see Remark\ \ref{rsecond}a)), and} \emph{%
therefore the arguments of the proof of \cite[Theorem 5.4]{rubsin} do not
work for the case of }$f^{c(\varphi )c(\varphi )^{\prime }}.$
\end{remark}

Another way of characterizing topical and anti-topical functions with the
aid of biconjugates is to combine conjugates with lower conjugates as
follows:

\begin{theorem}
\label{tantibiconj}If $(X,\mathcal{K})$\ is a pair satisfying $(A0^{\prime
}),(A1),$ then:

\emph{a) }A function $f:X\rightarrow \overline{\mathcal{K}}$ is topical if
and only if $f\not\equiv \top $ and 
\begin{equation}
f^{c(\varphi )\theta (\varphi )}=f.  \label{bic}
\end{equation}

\emph{b) }A function $f:X\rightarrow \overline{\mathcal{K}}$ is anti-topical
if and only if 
\begin{equation}
f^{\theta (\varphi )c(\varphi )}=f.  \label{antibic}
\end{equation}
\end{theorem}

\textbf{Proof}. a) If $f:X\rightarrow \overline{\mathcal{K}}$ is topical,
then $f\not\equiv \top $ (by (\ref{bbs-990})). Furthermore, $f^{c(\varphi )}$
is anti-topical (by Remark \ref{rsecond}a)), and hence by Theorem \ref
{tconj-antitop} applied to the function $f^{c(\varphi )},$ and Theorem \ref
{tconj}, we obtain 
\begin{equation*}
(f^{c(\varphi )})^{\theta (\varphi )}(y)=f^{c(\varphi
)}(y)^{-1}=(f(y)^{-1})^{-1}=f(y),\quad \quad \forall y\in X,
\end{equation*}
so $f$ satisfies (\ref{bic}).

Conversely, assume that $f$ satisfies $f\not\equiv \top $ and (\ref{bic}).
Then $f^{c(\varphi )}\not\equiv \varepsilon ,$ since otherwise by (\ref{bic}%
) and (\ref{counterp}) we would obtain $f(y)=$ $(f^{c(\varphi )})^{\theta
(\varphi )}(y)=\varepsilon ^{\theta (\varphi )}(y)=\top (y),\forall y\in X,$
in contradiction with the assumption $f\not\equiv \top .$ Consequently, by
Remark \ref{rsecond-bis}a), $f=(f^{c(\varphi )})^{\theta (\varphi )}$ is
topical.\quad \quad

b) If $f:X\rightarrow \overline{\mathcal{K}}$ is anti-topical, then $%
f\not\equiv \varepsilon $ (by (\ref{avem})). Consequently, by Remark \ref
{rsecond-bis}a), $f^{\theta (\varphi )}$ is topical and hence, by Theorem 
\ref{tconj} aplied to $f^{\theta (\varphi )},$ and Theorem \ref
{tconj-antitop} applied to $f,$ we obtain 
\begin{equation*}
f^{\theta (\varphi )c(\varphi )}(y)=f^{\theta (\varphi
)}(y)^{-1}=(f(y)^{-1})^{-1}=f(y),\quad \quad \forall y\in X.\quad \quad
\end{equation*}
\ \ 

Conversely, if $f$ satisfies (\ref{antibic}), then by Remark \ref{rsecond}
a), $f$ is anti-topical.\quad \quad $\square $\quad \quad\ 

\section{Polars of a set\label{s5}}

We shall study the following concept of ``polar set'' in our framework:

\begin{definition}
\label{dpolar}\emph{Let }$(X,\mathcal{K})$\emph{\ be a pair satisfying }$%
(A0^{\prime }),(A1),$\emph{\ let }$\mathcal{K}$ \emph{be enlarged to }$%
\overline{\mathcal{K}}:=\mathcal{K}\cup \{\top \}$\emph{\ as above, and let} 
\emph{\ }$\pi :X\times Y\rightarrow \overline{\mathcal{K}}$ \emph{be a
coupling function. The }$\pi $-polar \emph{of a subset }$G$\emph{\ of }$X$ 
\emph{is the subset }$G^{o(\pi )}$\emph{\ of} $Y$ \emph{defined by} 
\begin{equation}
G^{o(\pi )}:=\{y\in Y|\pi (g,y)\leq e,\forall g\in G\}=\{y\in Y|\sup_{g\in
G}\pi (g,y)\leq e\}.  \label{pi-polar}
\end{equation}
\emph{In particular, for the coupling function }$\pi =$\emph{\ }$\varphi
:X\times X\rightarrow \overline{\mathcal{K}}$\emph{\ of (\ref{coupl6}) we
have the }$\varphi $\emph{-polar of }$G$\emph{\ defined by} 
\begin{equation}
G^{o}:=G^{o(\varphi )}:=\{y\in Y|g/y\leq e,\forall g\in G\}=\{y\in
Y|\sup_{g\in G}(g/y)\leq e\},  \label{polarset}
\end{equation}
\emph{which we shall call simply the} polar \emph{of }$G.$
\end{definition}

\begin{remark}
\label{rR6.3}\emph{a) In particular, for }$Y=X,\pi :X\times Y\rightarrow 
\overline{\mathcal{K}}$ \emph{of (\ref{coupl6}) and any (non-empty) set }$%
G\subseteq X,$ \emph{by (\ref{polarset}) and (\ref{extenresid}) we have }$%
\inf X\notin G^{o}.$

\emph{b) Generalizing a concept of \cite{rubsin, radiant}, given two pairs }$%
(X,\mathcal{K}),(Y,\mathcal{K})$ \emph{satisfying assumptions }$(A0^{\prime
})$\emph{\ and }$(A1),$ \emph{a coupling function }$\pi :X\times
Y\rightarrow \overline{\mathcal{K}},$\emph{\ where} $\overline{\mathcal{K}}$%
\emph{\ is the minimal enlargement of }$\mathcal{K},$ \emph{and a subset} $G$
\emph{of} $X,$ \emph{we shall call} $\pi $-support function\emph{\ of }$G$ 
\emph{the function }$\sigma _{G,\pi }:Y\rightarrow \overline{\mathcal{K}}$%
\emph{\ defined by } 
\begin{equation}
\sigma _{G,\pi }(y):=\sup_{g\in G}\pi (g,y),\quad \quad \forall y\in Y.
\label{support}
\end{equation}
\emph{In the sequel we shall consider the particular case when }$Y=X$ \emph{%
and }$\pi =\varphi :X\times X\rightarrow \overline{\mathcal{K}}$\emph{\ of (%
\ref{coupl6}). In this case\ we shall call } 
\begin{equation}
\sigma _{G}(y):=\sigma _{G,\varphi }(y)=\sup_{g\in G}(g/y),\quad \quad
\forall y\in X,  \label{support3}
\end{equation}
\emph{the }support function\emph{\ of }$G.$ \emph{Note that for any set }$%
G\subseteq X,\sigma _{G}$ is an anti-topical function on\emph{\ }$X;$\emph{\
indeed, for each }$g\in G,$\emph{\ the function }$y\rightarrow g/y=\varphi
(g,y),\forall y\in X,$ \emph{is anti-topical (by Lemma \ref{lpartial}b)) and
the supremum of any family of anti-topical functions is anti-topical. Note
also that for }$y=\inf X$ \emph{we have} 
\begin{equation}
\sigma _{G}(\inf X)=\sup_{g\in G}(g/\inf X)=\top ,  \label{support4}
\end{equation}
\emph{and that the polar of} $G$ \emph{is the} $``e$\emph{-level set}$"$ 
\emph{of} $\sigma _{G},$ \emph{that is,} 
\begin{equation}
G^{o}=\{y\in X|\sigma _{G}(y)\leq e\}.  \label{lev}
\end{equation}

\emph{c) We recall that} \emph{a subset} $G$ \emph{of} $X$ \emph{is said to
be }downward\emph{, if} 
\begin{equation}
g\in G,x\in X,x\leq g\Rightarrow x\in G,  \label{downward}
\end{equation}
\emph{and }upward\emph{\ if} 
\begin{equation}
g\in G,x\in X,x\geq g\Rightarrow x\in G.  \label{upw}
\end{equation}

\emph{Let us observe that }for any set $G$ the polar $G^{o}=G^{o(\varphi )}$
is an upward set;\emph{\ indeed, if }$y\in G^{o},y^{\prime }\in X,y\leq
y^{\prime },$\emph{\ then }$y^{\prime }\in X\backslash \{\inf X\}$\emph{\
and }$g/y^{\prime }\leq g/y\leq e,$\emph{\ so }$y^{\prime }\in G^{o}.$\emph{%
\ Thus }$G^{o}$\emph{\ is upward. On the other hand, for the reflexion }$%
\overline{\varphi }$\emph{\ of} $\varphi $ \emph{(see (\ref{refl})), }the $%
\overline{\varphi }$-polar set 
\begin{equation}
G^{o(\overline{\varphi })}:=\{y\in X|y/g\leq e,\forall g\in G\}
\label{barpolar}
\end{equation}
is downward; \emph{indeed, if }$y\in G^{o(\overline{\varphi })},y^{\prime
}\leq y,$ \emph{then} $y^{\prime }/g\leq y/g\leq e,\forall g\in G.$
\end{remark}

We recall (see e.g. \cite{further} and the references therein) that for two
sets $X$ and $Y$ the\emph{\ dual\ }of any mapping\emph{\ }$\Delta
:2^{X}\rightarrow 2^{Y}$ (where $2^{X}$ denotes the family of all subsets of 
$X)$ is the mapping $\Delta ^{\prime }:2^{Y}\rightarrow 2^{X}$ defined by 
\begin{equation}
\Delta ^{\prime }(P):=\{x\in X|P\subseteq \Delta (\{x\})\},\quad \quad
\forall P\subseteq Y,  \label{dual2}
\end{equation}
Furthermore, we recall that a mapping $\Delta :2^{X}\rightarrow 2^{Y}$ is
called a \emph{polarity }(in \cite{further}, \cite{ACA} and earlier
references we have used the term \emph{duality}) if 
\begin{equation}
\Delta (G)=\cap _{g\in G}\Delta (\{g\}),\quad \quad \forall G\subseteq X.
\label{polarity}
\end{equation}
It is well-known and immediate that the dual $\Delta ^{\prime }$ of any
mapping $\Delta :2^{X}\rightarrow 2^{Y}$ is a polarity.

\begin{definition}
\label{dassocpol}\emph{If }$X,Y$\emph{\ are two sets and }$\pi :X\times
Y\rightarrow \overline{\mathcal{K}}$\emph{\ is a coupling function, the }%
polarity associated to $\pi ,$ \emph{or briefly, the} $\pi $-polarity, \emph{%
is the mapping }$\Delta _{\pi }:2^{X}\rightarrow 2^{Y}$\emph{\ defined by} 
\begin{equation}
\Delta _{\pi }(G):=G^{o(\pi )},\quad \quad \forall G\subseteq X;
\label{deltapi}
\end{equation}
\emph{with }$G^{o(\pi )}$ \emph{of (\ref{pi-polar}); the mapping }$\Delta
=\Delta _{\pi }$\emph{\ satisfies (\ref{polarity}), so it is indeed a
polarity.}
\end{definition}

We have the following extension of \cite[Lemma 2.1]{further}:

\begin{lemma}
\label{lL2.1}For the dual $\Delta _{\pi }^{\prime }$ of the mapping $\Delta
=\Delta _{\pi }\,$we have 
\begin{equation}
\Delta _{\pi }^{\prime }=\Delta _{\overline{\pi }},  \label{dual4}
\end{equation}
where $\overline{\pi }$ is the reflexion \emph{(\ref{refl})} of $\pi .$
\end{lemma}

\textbf{Proof}. By (\ref{dual2}), (\ref{deltapi}), (\ref{pi-polar}), (\ref
{refl}) and (\ref{deltapi}) (applied to $\overline{\pi }$ and $P)$, for any
subset $P$ of $Y$ we have 
\begin{eqnarray*}
\Delta _{\pi }^{\prime }(P) &=&\{x\in X|P\subseteq \Delta _{\pi
}(\{x\})\}=\{x\in X|\pi (x,y)\leq e,\forall y\in P\} \\
&=&\{x\in X|\overline{\pi }(y,x)\leq e,\forall y\in P\}=\Delta _{\overline{%
\pi }}(P).\quad \quad \square
\end{eqnarray*}

In the particular case when $\pi $ is a symmetric coupling function, that
is, satisfies (\ref{symm}), Lemma \ref{lL2.1} reduces to the following:

\begin{corollary}
\label{c2.1}If $(X,\mathcal{K})$ is a pair satisfying $(A0^{\prime }),(A1),$%
\ and $\pi :X\times X\rightarrow \overline{\mathcal{K}}$\ is a symmetric
coupling function, then $\Delta _{\pi }:2^{X}\rightarrow 2^{X}$ is
self-dual, that is, 
\begin{equation}
\Delta _{\pi }^{\prime }=\Delta _{\pi }.  \label{symm3}
\end{equation}
\end{corollary}

\begin{remark}
\emph{In particular, for }$\mathcal{K}=R_{\max }$ \emph{and }$X$\emph{\ an
arbitrary set, Corollary \ref{c2.1} has been obtained in \cite[Lemma 2.1]
{further}.}
\end{remark}

If $\Delta :2^{X}\rightarrow 2^{Y}$ is a polarity, with dual $\Delta
^{\prime },$ then by definition 
\begin{equation}
\Delta ^{\prime }\Delta (G):=\Delta ^{\prime }(\Delta (G)),\quad \quad
\forall G\subseteq X,  \label{bidelta}
\end{equation}
and a subset $G$ of $X$ is called $\Delta ^{\prime }\Delta $-\emph{convex }%
if $G=\Delta ^{\prime }\Delta (G),$ or equivalently (see \cite{ACA}), if for
each $x\notin G$ there exists $y\in Y$ such that 
\begin{equation}
G\subseteq \Delta ^{\prime }(\{y\}),\quad x\in X\backslash \Delta ^{\prime
}(\{y\}).  \label{deltaconv}
\end{equation}
In other words, a subset $G$ of $X$ is $\Delta ^{\prime }\Delta $-convex%
\emph{\ }if and only if it is $\mathcal{M}$-convex, i.e., ``convex with
respect to the family of sets'' $\mathcal{M}:=\{\Delta ^{\prime
}(\{y\})|y\in Y\}\subseteq X$ (see e.g. \cite{ACA} and the references
therein).

For the $\pi $-polarity $\Delta _{\pi }:2^{X}\rightarrow 2^{Y}$ of (\ref
{deltapi}), $\Delta _{\pi }^{\prime }\Delta _{\pi }$-convexity of a set $%
G\subseteq X$ means that for each $x\notin G$ there should exist $y\in Y$
such that $y\in \Delta _{\pi }(G)\backslash \Delta _{\pi }(\{x\}),$ or in
other words such that 
\begin{equation}
\sup_{g\in G}\pi (g,y)\leq e,\quad \pi (x,y)\nleq e.  \label{separ}
\end{equation}
$\newline
$For $Y=X$ and the coupling function $\pi =\varphi $ of (\ref{coupl6}), $%
\Delta _{\varphi }^{\prime }\Delta _{\varphi }$-convexity of a set $%
G\subseteq X$ means that for each $x\notin G$ there should exist $y\in X,$
or what is equivalent (by Remark \ref{rR6.3}a)), $y\in X\backslash \{\inf
X\},$\textbf{\ }such that 
\begin{equation}
\sup_{g\in G}(g/y)\leq e,\quad x/y\nleq e.  \label{separ2}
\end{equation}

In the particular case when $\leq $ is a total order on $\mathcal{K},$ (\ref
{separ}) and (\ref{separ2}) become, respectively, 
\begin{equation}
\sup_{g\in G}\pi (g,y)\leq e<\pi (x,y),\quad  \label{separ01}
\end{equation}
and 
\begin{equation}
\sup_{g\in G}(g/y)\leq e<x/y.  \label{separ21}
\end{equation}

Let us consider now the counterpart for ``bipolars'' $\Delta _{\varphi
}^{\prime }\Delta _{\varphi }(G)$ of the problem mentioned in Remark \ref
{rbiconj} (namely, to study the functions $f:X\rightarrow \overline{\mathcal{%
K}}$ such that $f^{c(\varphi )c(\varphi )^{\prime }}=f).$ Some results on $%
\Delta _{\varphi }^{\prime }\Delta _{\varphi }$-convex sets $G\subseteq X,$
without using this terminology, have been given in \cite{elem}. For example,
as in \cite{elem}, let us recall the following conditions on the relations
between the topologies of $\mathcal{K}$ and $X$ (for various other possible
conditions see \cite{AS} and the references therein):

(A2) For each $x\in X$ the function $u_{x}:\lambda \in \mathcal{K}%
\rightarrow \lambda x\in X$ is continuous, that is, for any $x\in X,\lambda
\in \mathcal{K}$ and any net $\{\lambda _{\kappa }\}\subset \mathcal{K}$
such that $\lambda _{\kappa }\rightarrow \lambda $ we have $\lambda _{\kappa
}x\rightarrow \lambda x$ (here we denote convergence both in $\mathcal{K}$
and in $X$ by $\rightarrow ,$ which will lead to no confusion).

(A3) For each $y\in X\backslash \{\inf X\}$ the function $y^{\diamond
}(.)=./y:X\rightarrow \mathcal{K}$ is continuous.

If (A2) holds, then a set $G\subseteq X$\ is said to be \emph{closed along
rays} if for each $x\in X$ the set 
\begin{equation}
H_{x}:=\{\lambda \in \mathcal{K}|\lambda x\in G\}\text{ }  \label{noua}
\end{equation}
is closed in $\mathcal{K}$ (that is, for any $x\in X$\ and any net $%
\{\lambda _{\kappa }\}\subset \mathcal{K}$ with $\{\lambda _{\kappa
}x\}\subset G$\ such that $\lambda _{\kappa }\rightarrow \lambda \in 
\mathcal{K}$\ we have $\lambda x\in G).$

We can recall now the following result of \cite[Theorem 10]{elem}:

\begin{theorem}
Let $(X$,$\mathcal{K})$ be a pair satisfying $(A0^{\prime }),(A1),(A2)$ and $%
(A3).$ For a subset $G$ of $X$ let us consider the following statements:

\emph{1}$^{\circ }$. $G$ is a closed downward set.

\emph{2}$^{\circ }$. $G$ is closed along rays and downward.

\emph{3}$^{\circ }$. We have $\Delta _{\varphi }^{\prime }\Delta _{\varphi
}(G)=G,$ that is, $G$ is $\Delta _{\varphi }^{\prime }\Delta _{\varphi }$%
-convex, where $\varphi $ is the coupling function \emph{(\ref{coupl6}).}

Then we have the implications $1^{\circ }\Rightarrow 2^{\circ }\Rightarrow
3^{\circ }.$

If the canonical order $\leq $ on $\mathcal{K}$ is total and if \emph{(A3)}
holds, then the statements \emph{1}$^{\circ }$, \emph{2}$^{\circ }$, and 
\emph{3}$^{\circ }$ are equivalent to each other and to the following
statements:

\emph{4}$^{\circ }.$ For each $x\in (X\backslash G)\backslash \{\inf X\}$
there exists $y\in X\backslash \{\inf X\}$ satisfying \emph{(\ref{separ21})}.

\emph{5}$^{\circ }.$ For each $x\in (X\backslash G)\backslash \{\inf X\}$
there exists $y\in X\backslash \{\inf X\}$ such that 
\begin{equation}
\sup_{g\in G}g/y<x/y.  \label{ulter}
\end{equation}
\end{theorem}

In particular, if $X=\mathcal{K}^{n},$\ where $\mathcal{K}=R_{\max },$\ then
the canonical order $\leq $ and the lattice operations in $\mathcal{K}$\ and 
$X$ coincide with the usual order and lattice operations on $\mathcal{K}^{n}$
and conditions (A2) and (A3) are satisfied for the order topology (by \cite
{CGQS}, Corollary 2.11), which is known to coincide with the usual
topologies on $\mathcal{K}$\ and $X$ (see \cite{CGQS}, the last observation
before Proposition 2.9).

\section{Support set. Support set at a point\label{s6}}

For a pair $(X,\mathcal{K})$ satisfying $(A0^{\prime }),(A1)$, we shall
denote by $\widetilde{\mathcal{T}}$ the set of all\emph{\ elementary topical}
functions $\widetilde{t}=\widetilde{t}_{y}:X\rightarrow \mathcal{K}$, that
is, all functions of the form 
\begin{equation}
\widetilde{t}(x)=\widetilde{t}_{y}(x):=x/y,\quad \quad \forall x\in
X,\forall y\in X\backslash \{\inf X\};  \label{elem1}
\end{equation}
we use here the notation $\widetilde{t}_{y}$ in order to avoid confusion
with the functions $t_{y}:X\rightarrow \overline{\mathcal{K}}$ of (\ref{unu}%
), (\ref{doijuma}).

We recall that for a set $X,$ a coupling function $\pi :X\times X\rightarrow
R_{\max },$ respectively $\pi :X\times X\rightarrow A,$ where $A=(A,\oplus
,\otimes )$ is a conditionally complete lattice ordered group, and a
function $f:$ $X\rightarrow \overline{R},$ respectively $f:$ $X\rightarrow 
\overline{A},$ the $\pi $-\emph{support set} of $f$ is the subset of $X$
defined \cite{rubsin, radiant} by 
\begin{equation}
\text{Supp}(f,\pi ):=\{y\in X|\pi (x,y)\leq f(x),\forall x\in X\}.
\label{suppset0}
\end{equation}

More generally, for a pair $(X,\mathcal{K})$ satisfying $(A0^{\prime }),(A1)$
and the coupling function $\pi =\varphi :X\times X\rightarrow \overline{%
\mathcal{K}}$ of (\ref{coupl6}) this leads us to define the following
concept of ``support set'' that will be suitable in our framework:

\begin{definition}
\label{dsuppset}\emph{Let} \emph{\ }$(X,\mathcal{K})$ \emph{be} \emph{a pair}
\emph{satisfying }$(A0^{\prime }),(A1).$ \emph{For a function }$%
f:X\rightarrow \overline{\mathcal{K}}$\emph{\ we shall call }support set%
\emph{\ }of\emph{\ }$f$\emph{\ (with respect to }$\widetilde{\mathcal{T}}$ $%
\ $\emph{of (\ref{elem1}}$))$\emph{\ the subset of }$X\backslash \{\inf X\}$ 
\emph{defined by} 
\begin{equation}
\text{Supp}(f,\widetilde{\mathcal{T}}):=\{y\in X\backslash \{\inf X\}|%
\widetilde{t}_{y}\leq f\}=\{y\in X\backslash \{\inf X\}|x/y\leq f(x),\forall
x\in X\}.  \label{suppset}
\end{equation}
\end{definition}

\begin{remark}
\label{r7.1}\emph{a)} For any function $f:X\rightarrow \overline{\mathcal{K}}%
,$ the support set Supp$(f,\widetilde{\mathcal{T}}\mathcal{)}$\emph{\ }is
upward. \emph{Indeed, if }$y\in $Supp$(f,\widetilde{\mathcal{T}}),y\leq
y^{\prime },$ \emph{then }$y^{\prime }\in X\backslash \{\inf X\}$ \emph{and }
\begin{equation*}
\widetilde{t}_{y^{\prime }}(x)=x/y^{\prime }\leq x/y=\widetilde{t}%
_{y}(x)\leq f(x),\quad \quad \forall x\in X,
\end{equation*}
\emph{\ so} $y^{\prime }\in $Supp$(f,\widetilde{\mathcal{T}}\mathcal{)}$.

\emph{b) If a function }$f:X\rightarrow \overline{\mathcal{K}}$\emph{\ is
topical then} 
\begin{equation}
f(x)=\max_{y\in \text{Supp}(f,\widetilde{\mathcal{T}}\mathcal{)}}\widetilde{t%
}_{y}(x),\quad \quad \forall x\in X.  \label{maxi}
\end{equation}
\emph{Indeed, by (\ref{suppset}) we have the inequality }$\geq $\emph{\ in (%
\ref{maxi}); furthermore, for any }$x\in X,$ \emph{the equality in (\ref
{maxi}) is attained for }$y=f(x)^{-1}x\in X\backslash f^{-1}(\varepsilon )$ 
\emph{when} $f(x)\in \mathcal{K}\backslash \{\varepsilon \},$ \emph{because
then} $x\in X\backslash \{\inf X\}$\emph{\ (by the topicality of }$f$ \emph{%
and} \emph{(\ref{bbs-990})) and } 
\begin{equation*}
\widetilde{t}_{y}(x)=\widetilde{t}_{f(x)^{-1}x}(x)=x/f(x)^{-1}x=f(x)x/x=f(x),
\end{equation*}
\emph{and for any }$y\in $Supp$(f,\widetilde{\mathcal{T}}\mathcal{)}$\emph{\
when} $f(x)=\varepsilon ,$ \emph{because then} $\widetilde{t}_{y}(x)=x/y\leq
f(x)=\varepsilon $ \emph{implies }$\widetilde{t}_{y}(x)=\varepsilon =f(x).$

\emph{From (\ref{maxi}) it follows that} for two topical functions $%
f_{1},f_{2}:X\rightarrow \overline{\mathcal{K}}$ we have 
\begin{equation}
f_{1}=f_{2}\Leftrightarrow \text{Supp}(f_{1},\widetilde{\mathcal{T}})%
\mathcal{=}\text{Supp}(f_{2},\widetilde{\mathcal{T}}\mathcal{)},
\label{univoc}
\end{equation}
\emph{which shows that }for a topical function $f:X\rightarrow \mathcal{K}$
the support set Supp$(f,\widetilde{\mathcal{T}})$\ determines uniquely $f.$%
\emph{\ This should be compared with \cite[Proposition 7.1]{radiant}.}
\end{remark}

Generalizing another concept studied in \cite{rubsin,radiant}, we introduce:

\begin{definition}
\label{datapoint}\emph{Let} \emph{\ }$(X,\mathcal{K})$ \emph{be} \emph{a pair%
} \emph{satisfying }$(A0^{\prime }),(A1),\pi :X\times X\rightarrow \overline{%
\mathcal{K}}$\emph{\ a coupling function and }$x_{0}\in X.$ \emph{For a
function }$f:X\rightarrow \overline{\mathcal{K}}$\emph{\ we shall call }%
support set\emph{\ }of\emph{\ }$f$ at $x_{0}$\emph{\ (with respect to }$%
\widetilde{\mathcal{T}}$ $\ $\emph{of (\ref{elem1}}$))$\emph{\ the subset of 
}$X\backslash \{\inf X\}$ \emph{defined by} 
\begin{equation}
\text{Supp}(f,\widetilde{\mathcal{T}};x_{0}):=\{y\in \text{Supp}(f,%
\widetilde{\mathcal{T}})|\pi (x_{0},y)=f(x_{0})\},  \label{atapoint1}
\end{equation}
\emph{that is,} 
\begin{equation}
\text{Supp}(f,\widetilde{\mathcal{T}};x_{0})=\{y\in X\backslash \{\inf
X\}|\pi (x,y)\leq f(x),\forall x\in X,\pi (x_{0},y)=f(x_{0})\}.
\label{atapoint0}
\end{equation}
\end{definition}

For a pair $(X,\mathcal{K})$ satisfying $(A0^{\prime }),(A1)$ and the
coupling functions $\pi =\varphi :X\times X\rightarrow \overline{\mathcal{K}}
$ and $\psi :X\times (X\times \overline{\mathcal{K}}\mathcal{)}\rightarrow 
\overline{\mathcal{K}}$ of (\ref{coupl6}) and (\ref{coupl}) respectively,
this leads us to define two concepts of ``support set at a point'' that will
be suitable in our framework: Supp$_{X}(f,\widetilde{\mathcal{T}};x_{0})$
and Supp$_{(X,\mathcal{K})}(f,\widetilde{\mathcal{T}};x_{0}).$

\begin{definition}
\label{dfirst}\emph{Let }$f:X\rightarrow \overline{\mathcal{K}}$\emph{\ be a
function and let }$x_{0}\in X.$\emph{\ We shall call }$X$-support set of $f$
at $x_{0}$ $\emph{the}$ $\emph{set}$\emph{\ } 
\begin{equation}
\text{Supp}_{X}(f,\widetilde{\mathcal{T}};x_{0}):=\{y\in X\backslash \{\inf
X\}|\widetilde{t}_{y}\leq f,\widetilde{t}_{y}(x_{0})=f(x_{0})\},
\label{subd0}
\end{equation}
\emph{where} $\widetilde{t}_{y}:X\rightarrow \overline{\mathcal{K}}$\emph{\
are the elementary topical functions (\ref{elem1}). In other words:} 
\begin{equation}
\text{Supp}_{X}(f,\widetilde{\mathcal{T}};x_{0})=\{y\in X\backslash \{\inf
X\}|x/y\leq f(x),\forall x\in X,x_{0}/y=f(x_{0})\}.  \label{000}
\end{equation}
\end{definition}

\begin{proposition}
\label{panter}Let $f:X\rightarrow \overline{\mathcal{K}}$ be a topical
function and let $x_{0}\in X$ be such that $f(x_{0})\in \mathcal{K}%
\backslash \{\varepsilon \}.$ Then 
\begin{equation}
x_{0}f(x_{0})^{-1}\in \text{Supp}_{X}(f,\widetilde{\mathcal{T}}%
;x_{0})\not=\emptyset .  \label{abram0}
\end{equation}
\end{proposition}

\textbf{Proof.} Since $f(x_{0})\in \mathcal{K}\backslash \{\varepsilon \},$
by (\ref{resid3/2}), the topicality of $f$ and Lemma \ref{lanti}a) we have 
\begin{equation}
x/f(x_{0})^{-1}x_{0}=f(x_{0})x/x_{0}\leq f(x),\quad \quad \forall x\in
X.\quad \quad  \label{abram}
\end{equation}

Hence defining 
\begin{equation}
y_{0}:=f(x_{0})^{-1}x_{0},  \label{abram2}
\end{equation}
we obtain $f(y_{0})=f(f(x_{0})^{-1}x_{0}))=e\neq \varepsilon ,$ so $y_{0}\in
X\backslash f^{-1}(\varepsilon )\subseteq X\backslash \{\inf X\}$ (by (\ref
{bbs-990})) and 
\begin{equation*}
x/y_{0}\leq f(x),\forall x\in X,\quad
x_{0}/y_{0}=x_{0}/f(x_{0})^{-1}x_{0}=f(x_{0}),
\end{equation*}
that is, $y_{0}\in $Supp$_{X}(f,\widetilde{\mathcal{T}};x_{0})\neq \emptyset 
$.\quad \quad $\square $

\begin{theorem}
\label{tprima}Let $f:X\rightarrow \overline{\mathcal{K}}$ be a topical
function and let $x_{0}\in X$ be such that $f(x_{0})\in \mathcal{K}%
\backslash \{\varepsilon \}.$ For an element $y\in X\backslash \{\inf X\}$
the following statements are equivalent:

\emph{1}$^{\circ }.$ $y\in $Supp$_{X}(f,\widetilde{\mathcal{T}};x_{0}).$

\emph{2}$^{\circ }.$ We have 
\begin{equation}
f(y)=e,\widetilde{t}_{y}(x_{0})=f(x_{0}).  \label{subd1}
\end{equation}

\emph{3}$^{\circ }.$ We have 
\begin{equation}
f(y)=e,\widetilde{t}_{y}(x_{0})\geq f(x_{0}).  \label{main-eq-200}
\end{equation}
\end{theorem}

\textbf{Proof.} 1$^{\circ }\Rightarrow 2^{\circ }.$ Assume that $y\in $Supp$%
_{X}(f,\widetilde{\mathcal{T}};x_{0})$. Then by (\ref{subd0}) and (\ref
{resid5}) we have 
\begin{equation}
yf(x_{0})=y\widetilde{t}_{y}(x_{0})=y(x_{0}/y)\leq x_{0}.  \label{sem-eq-33}
\end{equation}

Hence, since $f$ is topical, we obtain 
\begin{equation}
f(y)f(x_{0})=f(yf(x_{0}))\leq f(x_{0}).  \label{eq-202}
\end{equation}
But, by our assumption we have $f(x_{0})\in \mathcal{K}\backslash
\{\varepsilon \},$ whence multiplying (\ref{eq-202}) with $f(x_{0})^{-1},$
we obtain 
\begin{equation}
f(y)\leq e.  \label{bfirst-ineq-1}
\end{equation}

On the other hand, using the fact that by (\ref{resid4}) we have 
\begin{equation*}
\widetilde{t}_{y}(y)=y/y=e,\quad \quad \forall y\in X\backslash \{\inf X\},
\end{equation*}
and again that $y\in $Supp$_{X}(f,\widetilde{\mathcal{T}};x_{0})$, one has: 
\begin{equation}
e=\widetilde{t}_{y}(y)\leq f(y).  \label{ineq-11}
\end{equation}

Combining (\ref{bfirst-ineq-1}) and (\ref{ineq-11}) yields $f(y)=e.$

Finally, from $y\in $Supp$_{X}(f,\widetilde{\mathcal{T}};x_{0})$ and $f(y)=e$
it follows that we have 2$^{\circ }.$

The implication 2$^{\circ }\Rightarrow 3^{\circ }$ is obvious.

3$^{\circ }\Rightarrow 1^{\circ }.$ Assume that $y\in X\backslash \{\inf X\}$
satisfies 3$^{\circ }.$ Then using $f(y)=e$, the topicality of $f$ and Lemma 
\ref{lanti}a), one has 
\begin{equation}
\widetilde{t}_{y}(x)=e(x/y)=f(y)(x/y)\leq f(x),\quad \quad x\in X.
\label{bff-01}
\end{equation}
\ \ \ \ \ \ \ \ \ \ \ \ \ \ \ \ \ \ \ \ \ \ \ \ \ \ \ \ \ \ \ \ \ \ \ \ \ \
\ \ \ \ \ \ \ \ \ \ \ \ \ \ \ \ \ \ \ \ \ \ \ \ \ \ \ \ \ \ \ \ \ \ \ \ \ \
\ \ \ \ \ \ \ \ \ \ \ \ \ \ \ \ \ \ \ \ \ \ \ \ \ \ \ \ \ \ \ \ \ \ \ \ \ \
\ \ \ \ \ \ \ \ \ \ \ \ \ \ \ \ \ \ \ \ \ \ \ \ \ \ \ \ \ \ \ \ \ \ \ \ \ \
\ \ \ \ \ \ \ \ \ \ \ \ \ \ \ \ \ \ \ \ \ \ \ \ \ \ \ \ \ \ \ \ \ \ \ \ \ \
\ \ \ \ \ \ \ \ \ \ \ \ \ \ \ \ \ \ \ \ \ \ \ \ \ 

On the other hand, by $3^{\circ }$ we have 
\begin{equation}
\widetilde{t}_{y}(x_{0})=x_{0}/y\geq f(x_{0}),  \label{bff-03}
\end{equation}
and hence, by (\ref{bff-01}) (for $x=x_{0})$ and (\ref{bff-03}), we obtain 
\begin{equation}
\widetilde{t}_{y}(x_{0})=x_{0}/y=f(x_{0}).  \label{bff-04}
\end{equation}

From (\ref{bff-04}) and (\ref{bff-01}) it follows that $y\in $Supp$_{X}(f,%
\widetilde{\mathcal{T}};x_{0})$.\quad \quad $\square $

\begin{definition}
\label{dsecond}\emph{Let }$f:X\rightarrow \overline{\mathcal{K}}$\emph{\ be
a function and let }$x_{0}\in X.$\emph{\ We shall call }$(X,\mathcal{K})$%
-support set of $f$ at $x_{0}$ $\emph{the}$ $\emph{set}$ 
\begin{equation}
\text{Supp}_{(X,\mathcal{K})}(f,\widetilde{\mathcal{T}};x_{0}):=\{(y,d)\in
(X\backslash \{\inf X\})\times \mathcal{K}|s_{y,d}\leq
f,s_{y,d}(x_{0})=f(x_{0})\},  \label{subd}
\end{equation}
\emph{where} $s_{y,d}:X\rightarrow \overline{\mathcal{K}}$ \emph{are the
functions (\ref{aff1}),(\ref{aff7}),(\ref{aff5})}.
\end{definition}

\begin{proposition}
\label{pdupa}Let $f:X\rightarrow \overline{\mathcal{K}}$ be a topical
function and let $x_{0}\in X$ be such that $f(x_{0})\in \mathcal{K}%
\backslash \{\varepsilon \}.$ Then 
\begin{equation}
(f(x_{0})^{-1}x_{0},f(x_{0}))\in \text{Supp}_{(X,\mathcal{K})}(f,\widetilde{%
\mathcal{T}};x_{0})\neq \emptyset .  \label{abram6}
\end{equation}
\end{proposition}

\textbf{Proof}. If $f$ is topical and $f(x_{0})\in \mathcal{K}\backslash
\{\varepsilon \},$ then $(f(x_{0})^{-1}x_{0},f(x_{0}))\in (X\backslash
\{\inf X\})\times (\mathcal{K}\backslash \{\varepsilon \})$ and by (\ref
{resid3/2}), the topicality of $f$ and Theorem \ref{lanti-bis1}a) we have

$_{{}}$%
\begin{eqnarray}
\inf \{x/f(x_{0})^{-1}x_{0},d\} &=&\inf \{f(x_{0})x/x_{0},d\}\leq \inf
\{f(x),d\}  \notag \\
&\leq &f(x),\quad \quad \forall x\in X,\forall d\in \mathcal{K}.\quad \quad
\label{abram3}
\end{eqnarray}
Consequently, 
\begin{eqnarray*}
s_{f(x_{0})^{-1}x_{0},f(x_{0})}(x) &=&\inf
\{x/f(x_{0})^{-1}x_{0},f(x_{0})\}\leq f(x),\quad \quad \forall x\in X, \\
s_{f(x_{0})^{-1}x_{0},f(x_{0})}(x_{0}) &=&\inf
\{x_{0}/f(x_{0})^{-1}x_{0},f(x_{0})\}=f(x_{0}),
\end{eqnarray*}
that is, $(f(x_{0})^{-1}x_{0},f(x_{0}))\in $ Supp$_{(X,\mathcal{K})}(f,%
\widetilde{\mathcal{T}};x_{0})\neq \emptyset .\quad \quad \square $

\begin{theorem}
\label{tadoua}Let $f:X\rightarrow \overline{\mathcal{K}}$ be a topical
function and let $x_{0}\in X$ be such that $f(x_{0})\in \mathcal{K}%
\backslash \{\varepsilon \}.$ For a pair $(y,d)\in (X\backslash \{\inf
X\})\times \mathcal{K},$ the following statements are equivalent:

\emph{1}$^{\circ }.$ $(y,d)\in $Supp$_{(X,\mathcal{K})}(f,\widetilde{%
\mathcal{T}};x_{0}).$

\emph{2}$^{\circ }.$ We have 
\begin{equation}
f(y)=e,s_{y,d}(x_{0})=f(x_{0}).  \label{subd2}
\end{equation}

\emph{3}$^{\circ }.$ We have 
\begin{equation}
f(y)=e,s_{y,d}(x_{0})\geq f(x_{0}).  \label{main-eq-100}
\end{equation}
\end{theorem}

\textbf{Proof.} 1$^{\circ }\Rightarrow 2^{\circ }.$ Assume that $(y,d)\in $
Supp$_{(X,\mathcal{K})}(f,\widetilde{\mathcal{T}};x_{0})$. Then by (\ref
{subd}) we have 
\begin{equation}
f(x_{0})=s_{y,d}(x_{0})=\inf \{x_{0}/y,d\}\leq x_{0}/y,  \label{main-eq-101}
\end{equation}
whence by (\ref{resid5}), 
\begin{equation*}
yf(x_{0})\leq y(x_{0}/y)\leq x_{0}.
\end{equation*}
Hence, since $f$ is topical, we obtain

\begin{equation}
f(y)f(x_{0})=f(yf(x_{0}))\leq f(x_{0}).  \label{eq-102}
\end{equation}
But, by our assumption we have $f(x_{0})\in \mathcal{K}\backslash
\{\varepsilon \},$ whence multiplying (\ref{eq-102}) with $f(x_{0})^{-1}$ we
obtain

\begin{equation}
f(y)\leq e.  \label{ineq-1}
\end{equation}

On the other hand, using again that $(y,d)\in $ Supp$_{(X,\mathcal{K})}(f,%
\widetilde{\mathcal{T}};x_{0}),$ we have 
\begin{equation*}
\inf \{x/y,d\}=s_{y,d}(x)\leq f(x),\quad \quad \forall x\in X,
\end{equation*}
which for $x=dy$ gives, using (\ref{resid4}) and that $f$ is homogeneous, 
\begin{equation*}
d=\inf \{dy/y,d\}\leq f(dy)=df(y).
\end{equation*}
Hence, multiplying with $d^{-1}$ (which exists by $\varepsilon
<f(x_{0})=\inf \{x_{0}/y,d\}\leq d$ ), we get 
\begin{equation}
e\leq f(y).  \label{ineq-2}
\end{equation}

Combining (\ref{ineq-1}) and (\ref{ineq-2}) yields $f(y)=e.$

Finally, from $(y,d)\in $ Supp$_{(X,\mathcal{K})}(f,\widetilde{\mathcal{T}}%
;x_{0})$ and $f(y)=e$ it follows that we have 2$^{\circ }$.

The implication 2$^{\circ }\Rightarrow 3^{\circ }$ is obvious.

3$^{\circ }\Rightarrow 1^{\circ }.$ Assume that the pair $(y,d)\in
(X\backslash \{\inf X\})\times \mathcal{K}$ satisfies 3$^{\circ }.$ Then
using $f(y)=e$, the topicality of $f$ and Theorem \ref{lanti-bis1}a), we
obtain 
\begin{equation}
s_{y,d}(x)=es_{y,d}(x)=f(y)s_{y,d}(x)\leq f(x),\quad \quad \forall x\in X.
\label{cl}
\end{equation}

On the other hand, by 3$^{\circ }$ we have $s_{y,d}(x_{0})\geq f(x_{0}),$
and hence the equality $s_{y,d}(x_{0})=f(x_{0}),$ which, by (\ref{cl}), \
yields $(y,d)\in $ Supp$_{(X,\mathcal{K})}(f,\widetilde{\mathcal{T}}%
;x_{0}).\quad \square $

\begin{remark}
\label{rsubdiff}\emph{We recall that for a set} $X,$\emph{\ a coupling
function }$\pi :X\times X\rightarrow A,$ \emph{where} $A=(A,\oplus ,\otimes
) $ \emph{is a conditionally complete lattice ordered group, a function} $f:$
$X\rightarrow \overline{A}$ \emph{(where }$\overline{A}$ \emph{is the
minimal enlargement of} $A)$ \emph{and} \emph{\ a point }$x_{0}\in X$ \emph{%
with} $f(x_{0})\in A,$\emph{\ the} $\pi $-subdifferential of $f$ at $x_{0}$ 
\emph{is the subset of }$X$ \emph{defined \cite{radiant,MS} by} 
\begin{equation}
\partial _{\pi }f(x_{0}):=\{y_{0}\in X|\pi (x,y_{0})\otimes \pi
(x_{0},y_{0})^{-1}\otimes f(x_{0})\leq f(x),\forall x\in X\}.
\label{pisubd2}
\end{equation}

\emph{For} $\pi =\pi _{\mu }$ \emph{of (\ref{scal}) some properties of }$%
\partial _{\pi _{\mu }}f(x_{0}),$ \emph{e.g. that it contains the }$\pi $%
\emph{-support set of }$f$\emph{\ at} $x_{0}$\emph{\ } \emph{and some other
connections between these sets} \emph{have been proved in} \emph{
\cite[Section 8]{radiant}. This suggests to attempt to introduce,} \emph{in
our framework of a pair} $(X,\mathcal{K})$ \emph{satisfying} $(A0^{\prime
}),(A1)$ \emph{and the coupling function} $\varphi :X\times X\rightarrow 
\overline{\mathcal{K}}$ \emph{of} \emph{(\ref{coupl6}), the following
concept: We shall call }$\varphi $-subdifferential of a function $%
f:X\rightarrow \overline{\mathcal{K}}$ at a point $x_{0}$ with $f(x_{0})\in 
\mathcal{K}$, \emph{the subset of} $X$ \emph{defined by} 
\begin{equation}
\partial _{\mathcal{\varphi }}f(x_{0}):=\{y_{0}\in
X|(x/y_{0})(x_{0}/y_{0})^{-1}f(x_{0})\leq f(x),\forall x\in X\}.
\label{pisubd3}
\end{equation}
\emph{The proofs of the properties of }$\partial _{\pi _{\mu }}f(x_{0})$ 
\emph{given in \cite[Section 8]{radiant} lean heavily on the fact that for
each }$x\in X$ \emph{the partial function }$\pi _{\mu }(x,.)$ \emph{is} 
\emph{topical, but in the case of the coupling function} $\varphi :X\times
X\rightarrow \overline{\mathcal{K}}$ \emph{of (\ref{coupl6}) for each }$x\in
X$ \emph{the partial function }$\varphi (x,.)$ \emph{is anti-topical (see
Lemma \ref{lpartial}b)). Therefore those proofs of \cite{radiant} cannot be
adapted to our general framework.}
\end{remark}

\textbf{Acknowledgement\bigskip }

Viorel Nitica was partially supported by a grant from Simons Foundation
208729.\medskip\ 

We would like to thank the anonymous referee for careful reading of this
paper and for making useful suggestions.

\end{document}